\documentstyle[a4,leqno,amsfonts,12pt]{article}

\title{CONSTRUCTIVE FUNCTION THEORY ON SETS OF THE COMPLEX
PLANE THROUGH POTENTIAL THEORY AND GEOMETRIC FUNCTION THEORY}
\def\shorttitle{Constructive function theory}

\author{V. V. Andrievskii}
\def\shortauthor{Andrievskii}

\def\versiondate{12/13/2005}
\def\startpagenumber{1}
\def\volumenumber{1}
\def\year{2005}
\newcommand{\beginddoc}{
\begin{document}
\maketitle \insert\footins{\scriptsize
\medskip
\baselineskip 8pt \leftline{Surveys in Approximation Theory}
\leftline{Volume \volumenumber, \year.
pp.~\thepage--\pageref{endpage}.} \leftline{Copyright \copyright\
2006 Surveys in Approximation Theory.} \leftline{ISSN 1555-578X}
\leftline{All rights of reproduction in any form reserved.}
\smallskip
\par\allowbreak}
\tableofcontents} \markboth{{\it \shortauthor}}{{\it \shorttitle}}
\markright{{\it \shorttitle}}
\def\endddoc{\label{endpage}\end{document}}
\date{{\small \versiondate}}


\newcommand{\C}{{\bf C}}
\newcommand{\ol}{\overline}
\newcommand{\OD}{\overline{\bf D}}
\newcommand{\la}{\lambda}
\def\MG{m_2(\overline{G})}
\newcommand{\Da}{{\bf D}}
\newcommand{\fk}{\mbox{\bf k}}
\newcommand{\kap}{\mbox{cap}}
\newcommand{\cL}{{\cal{E}}}
\newcommand{\rbox}{$\:\:$ \raisebox{-1ex}{$\:\Box\:$}}
\newcommand{\OC}{\overline{\bf C}}
\newcommand{\R}{{\bf R}}
\newcommand{\He}{{\bf H}}
\newcommand{\Pn}{{\bf P}_n}
\newcommand{\T}{{\bf T}}
\newcommand{\Pa}{{\bf P}}
\newcommand{\N}{{\bf N}}
\newcommand{\bP}{{\bf {\Pi}}}
\newcommand{\D}{{\bf D}}
\newcommand{\de}{\delta}
\newcommand{\De}{\Delta}
\newcommand{\mb}{\mbox}
\newcommand{\fn}{\mbox{\bf n}}
\newcommand{\beq}{\begin{equation}}
\newcommand{\eeq}{\end{equation}}
\newcommand{\oge}{\succeq}
\newcommand{\ole}{\preceq}
\newcommand{\ve}{\varepsilon}
\newcommand{\bl}{\backslash}
\newcommand{\ov}{\overline}
\newcommand{\al}{\alpha}
\newcommand{\Ze}{{\cal Z}}
\newcommand{\Z}{{\bf Z}}
\newcommand{\CC}{\overline{\bf C}}
\newcommand{\be}{\beta}
\newcommand{\Om}{\Omega}
\newcommand{\om}{\omega}
\newcommand{\z}{\zeta}
\newcommand{\ka}{\kappa}
\newcommand{\ga}{\gamma}
\newcommand{\Ga}{\Gamma}
\newcommand{\si}{\sigma}
\newtheorem{thm}{Theorem}[section]
\newtheorem{cor}[thm]{Corollary}
\newtheorem{lem}[thm]{Lemma}

\beginddoc


\begin{abstract}
This is a survey of some recent results
 concerning polynomial inequalities and polynomial approximation
 of functions in the complex plane.
 The
results are achieved  by the application of methods and techniques
of modern geometric function theory and potential theory.
\end{abstract}



\section{Introduction}

Constructive function theory, or more generally, the theory of the
representation of functions by series of polynomials and rational
functions, may be described as part of the intersection  of the
analysis and applied mathematics. The main feature  of the
research discussed in this survey concerns  new methods based on
conformal invariants to solve problems arising in potential
theory, geometric function theory and approximation theory.

The  harmonic measure, module and extremal length of a family of
curves   serve as the main tool. A significant part of the work
depends on new techniques for the study of the special conformal
mapping  of the upper half-plane onto the upper half-plane with
vertical slits. These techniques have independent value and have
already been  applied to other areas of mathematics.

This survey is organized as follows.
 Section 2 is devoted to the
properties of the Green function $g_{\CC\setminus E}$ and
equilibrium measure $\mu_E$ of a compact set $E$ on the real line
$\R$.
  Recently, Totik \cite{tot},  Carleson and Totik \cite{cartot},
  and the author \cite{andfni, andh, ands}
  suggested  new methods to approach these objects.
  We use a new representation of basic
notions of potential
 theory (logarithmic capacity, the Green function, and
 equilibrium measure) in terms of a conformal mapping of the
 exterior of the  interval $[0,1]$  onto the exterior of the unit
 disk $\D$
 with finite or infinite number of radial slits
 \cite{andup} - \cite{andh}. This method provides a number of new
 links between potential theory and the theory of univalent functions.
 Later in this section, we describe the connection between uniformly perfect compact
sets and John domains.
  We give a new interpretation (and a generalization)
of a recent remarkable result by  Totik \cite[(2.8) and
(2.12)]{tot} concerning the smoothness properties of $g_\Om$ and
$\mu_E$.  We also demonstrate that if for $E\subset[0,1]$ the
Green function satisfies the $1/2$-H\"older condition locally at
the origin, then the density of $E$ at $0$, in terms of
logarithmic capacity, is the same as that of  the whole interval
$[0,1]$. We analyze the geometry of Cantor-type sets and propose
an
 extension of the results by Totik
 \cite[Theorem 5.3]{tot} on Cantor-type sets possessing the $1/2$-H\"older
 continuous Green
 function. We also construct two examples of sets  of  minimum Hausdorff
dimension with  Green function  satisfying the $1/2$-H\"older
condition either uniformly or locally.

In  Section 3 we continue to discuss the properties of the Green
function,
 but now we motivate this investigation by deriving
 Remez-type polynomial inequalities.
 We give   sharp uniform bounds for exponentials of
logarithmic potentials if the logarithmic capacity of the subset,
where they are at most $1$, is known.
 We also
propose a technique to derive Remez-type inequalities for complex
polynomials.
 The known results in this direction are
scarce
 and they are proved for relatively
simple geometrical cases by using  methods of real analysis. We
propose to use modern methods  of complex analysis, such as  the
application of conformal invariants in constructive
 function theory   and the theory of quasiconformal mappings in the
plane, to study  metric properties of complex polynomials.
 Based on this idea we discuss a number of
problems motivated by \cite{erdlisaf}.

  In Section 4
 we  consider several applications of methods and techniques
covered in the previous two sections to questions arising in
constructive function theory. The main idea of this section is to
create a link between potential theory, geometric function theory
and approximation theory.
 We present a new necessary condition and a new sufficient condition for the
approximation of the reciprocal of an entire function by
reciprocals of polynomials on the non-negative real line with
geometric speed of convergence.
 The Nikol'skii-Timan-Dzjadyk theorem concerning
polynomial approximation of functions on the interval $[-1,1]$ is
generalized to the case of approximation of functions given on a
compact set on the real line.
 For analytic functions defined on a continuum $E$ in the complex
plane, we discuss Dzjadyk-type polynomial approximations   in
terms of the $k$-th modulus of continuity ($k\ge 1$) with
simultaneous interpolation at given points of $E$ and decaying
strictly inside as $e^{-cn^\alpha}$, where $c$ and $\alpha$ are
positive constants independent of the degree $n$ of the
approximating polynomial.

Each section concludes with a list of open problems.

\section{Potential theory}

\subsection{Basic conformal mapping}
 Let $E\subset \C$ be a compact
set of positive logarithmic capacity cap$(E)$ with connected
complement $\Om:=\OC\setminus E$
 with respect to the extended complex plane
$\OC=\C\cup\{\infty\}$,
 $g_\Om(z)=g_\Om(z,\infty)$ be the Green function of $\Om$ with pole
 at infinity, and $\mu_E$ be the equilibrium measure for the set $E$
 (see \cite{lan} and \cite{saftot} for further details on
 logarithmic potential theory). The metric properties of $g_\Omega$ and
 $\mu_E$ are  of independent interest in potential
 theory (see, for example, \cite{car,maz, lit, saftot, andbla, cartot, tot, andfni, andh}).
 They also
  play an important role in  problems concerning
 polynomial approximation of continuous functions on $E$ (see, for example,
  \cite{tam,dzj,gai,shi1,andbeldzj})
  and the behavior of
 polynomials with a known uniform norm along $E$ (see, for example,
  \cite{wal,plec,plem,biavol,erdlisaf,borerd,tot1,tot2}).

 Note that sets in $\R$ present an important special case of
 general sets in $\C$. This, for instance, is due to the following standard way to
 simplify problems concerning estimation of the Green function and
 capacity.
 For $E\subset \C$ denote by
 $E_*:=\{ r:\, \{|z|=r\}\cap E\neq\emptyset\}$
 the circular projection of $E$ onto the non-negative real line $\R^+:=\{x\in\R:\, x\ge0\}$. Then
 $$\mb{cap}(E)\ge \mb{cap}(E_*)$$
 and
 $$g_{\OC\setminus E}(-x)\le g_{\OC\setminus E_*}(-x),\quad x>0$$
 (provided that cap$(E_*)>0$). That is, among those sets that have a
 given circular projection $E_*\subset\R^+$ the smallest capacity
 occurs for $E=E_*$ and the worst behavior of the Green function occurs
 for the same $E=E_*$.

  In this survey we discuss a number of  problems in potential
   theory, polynomial inequalities, and constructive function theory for the case
 where $E$ is a subset of  $\R$.

  The main
 idea of our approach is to connect $g_\Om,\mu_E,$ and cap$(E)$
 with the special conformal mapping $F=F_E$ described below. This
 conformal mapping was recently investigated  in
  \cite{andup} - \cite{andh} (written in another form it was also discussed in
  \cite{wid,lev, sodyud}).

 Let $E\subset [0,1] $ be a regular set such that $0\in E,\, 1\in E$. Then
 $[0,1]\setminus E=\sum_{j=1}^N(a_j,b_j)$, where $N$ is finite or
 infinite.

 Denote  by $\He:=\{z:\, \Im(z)>0\}$ the
upper half-plane and consider the function \beq
F(z)=F_E(z):=\exp\left(\int_E\log(z-\z)\,
d\mu_E(\z)-\log\kap(E)\right), \quad z\in \He.\label{fdef}\eeq
 It is analytic in
$\He$.

 Since
  $$ g_{\Om}(z)=\log\frac{1}{\kap
(E)}-\int\log\frac{1}{|z-t|}d\mu_E(t),\quad z\in \Om,
 $$
 the function $F$ has the following obvious
properties: $$|F(z)|=e^{g_\Om(z)}>1,\quad z\in\He,$$
$$\Im(F(z))=e^{g_\Om(z)}\sin\left(\int_E\arg(z-\z)\,
d\mu_E(\z)\right)>0,\quad z\in\He.$$ Moreover, $F$ can be extended
from $\He$ continuously to $\ov{\He}$ such that
 \begin{eqnarray*}
 |F(z)|&=&1,\quad
z\in E,\\ F(x)&=&e^{g_\Om(x)}>1,\quad x>1,\\
F(x)&=&-e^{g_\Om(x)}<-1,\quad x<0.
 \end{eqnarray*}
  Next, for any
$1\le j\le N$ and $a_j\le x_1<x_2\le b_j$, we have
$$\arg\left(\frac{F(x_2)}{F(x_1)}\right)=
\arg\exp\left(\int_E\log\frac{x_2-\z}{x_1-\z}\,
d\mu_E(\z)\right)=0,$$ that is, $$\arg F(x_1)=\arg F(x_2),\quad
a_j\le x_1<x_2\le b_j.$$

Our next objective is to show that $F$ is univalent in $\He$. We
shall use the following simple result. Let $\sqrt{z^2-1},\,
z\in\OC\setminus [-1,1],$ be the analytic function defined in a
neighborhood of infinity as $$
\sqrt{z^2-1}=z\left(1-\frac{1}{2z^2}+\cdots\right).$$ Then
 for any $-1\le x\le 1$ and $z\in\He$,
\beq\label{dd}u_x(z):= \Re\left(\frac{\sqrt{z^2-1}}{z-x}\right)\ge
0.\eeq
 Using the reflection principle we can extend
$F$ to a function analytic in $\OC\setminus[0,1]$ by the formula
$$F(z):=\ov{F(\ov{z})}, \quad z\in\C\setminus \ov{\He},$$ and
consider the function $$h(w):= \frac{1}{F(J(w))}\, ,\quad
w\in\Da:=\{w:\,|w|<1\},$$ where $J$ is a linear transformation of
the Joukowski mapping, namely
$$J(w):=\frac{1}{2}\left(\frac{1}{2}\left(w+\frac{1}{w}\right)+1\right),$$
which maps the unit disk $\Da$ onto $\OC\setminus[0,1]$. Note that
the inverse mapping is defined as follows
$$w=J^{-1}(z)=(2z-1)-\sqrt{(2z-1)^2-1},\quad z\in
\OC\setminus[0,1].$$ Therefore, for $z\in\He$ and $w=J^{-1}(z)\in
\Da$, we obtain \begin{eqnarray*}\frac{wh'(w)}{h(w)}&=& w(\log
h(w))'=-w\left(\int_E\log(J(w)-\z)\, d\mu_E(\z)\right)'\\ &=& -w\,
J'(w)\int_E\frac{d\mu_E(\z)}{z-\z}=-\frac{1}{4}\left(w-\frac{1}{w}\right)
\int_E\frac{d\mu_E(\z)}{z-\z}\\
&=&\frac{1}{2}\int_E\frac{\sqrt{(2z-1)^2-1}}{z-\z}\,
d\mu_E(\z)=\int_E \frac{\sqrt{(2z-1)^2-1}}{(2z-1)-(2\z-1)}\,
d\mu_E(\z).\end{eqnarray*} According to (\ref{dd})
 for $w$
under consideration we have
$$\Re\left(\frac{wh'(w)}{h(w)}\right)\ge 0. $$ Because of the
symmetry and the maximum principle for harmonic functions we
obtain $$\Re\left(\frac{wh'(w)}{h(w)}\right)> 0,\quad w\in\Da. $$
It means that $h$ is a conformal mapping of $\Da$ onto a starlike
domain (cf. \cite[p. 42]{pom}).

Hence, $F$ is univalent and maps $\OC\setminus[0,1]$ onto a  (with
respect to $\infty$) starlike domain $\OC\setminus K_E$  with the
following properties: $\OC\setminus K_E$ is symmetric with respect
to the real line and coincides with the exterior of the unit disk
with $2N$ slits.

Note that
\begin{eqnarray}
 \kap(E)&=&\frac{1}{4\,\kap(K_E)},\nonumber\\
 \label{1111}
 g_\Om(z)&=&\log|F(z)|,\quad z\in
\Om,\\
 \pi\mu_E([a,b])&=&|F([a,b]\cap E)|,\nonumber
 \end{eqnarray}
 where $|A|$ denotes the linear Lebesgue measure (length) of a
 Borel set $A\subset\C$.

The connection between the geometry of $E$ and the properties of
the conformal mapping $F$   can be studied using
  conformal invariants such as the
 extremal length and module of a family of curves
  (see \cite{ah, lv, pom1}).

  Below, we describe some typical results of this investigation.

\subsection{Uniformly perfect subsets of the real line and John
domains}
 The uniformly perfect sets in the complex plane $\C$,
introduced by Beardon and Pommerenke \cite{pombea}, are defined as
follows. A compact set $E\subset\C$ is  {\it uniformly perfect} if
there exists a  constant $0<c<1$ such that for all $z\in E$:
 $$
 E\cap\{ \z:\, cr\le|z-\z|\le r\}\neq\emptyset,\quad
0<r<\mb{diam}(E):=\sup_{z,\z\in E}|z-\z|.
 $$
 Uniformly perfect sets arise in many areas of complex analysis.
 For example, many results for simply connected domains can be
 extended to domains with uniformly perfect boundary (see, for example,
 \cite{pomu,pom5,zhe,bis}).
Pommerenke \cite{pomu} has shown that uniformly perfect sets can
be described using a density condition expressed in terms of the
logarithmic capacity. Namely, $E$ is uniformly perfect if and only
if there exists a positive constant $c$ such that for all $z\in
E$:
 \beq\label{2.21.2} \kap(E\cap\{\z:\, |\z-z|\le r\})\ge c\, r,\quad 0<r\le\mb{ diam}(E).
 \eeq
 It
follows immediately from (\ref{2.21.2}) that each component of
$\OC\setminus E$ is regular (for the Dirichlet problem).

Note that sets $E$ with connected complement $\OC\setminus E$
satisfying (\ref{2.21.2}) play a significant role in the solution
of the inverse problem of the constructive theory of functions of
a complex variable. We refer to \cite{tam} where they are called
{\it $c$-dense} sets.

Another remarkable geometric condition used in direct theorems of
approximation theory in $\C$ (cf. \cite{gai, and93, bagboslev})
defines  a John domain \cite{marsar, pom1}. We consider only the
case of  a simply connected domain $\Om\subset\OC$ such that
$\infty\in\Om$. Following \cite[p. 96]{pom1} we call $\Om$ a {\it
John domain} if there exists a positive constant $c$ such that for
every rectilinear crosscut $[a,b]$ of $\Om$, $$ \mb{diam} (H)\le
c|a-b| $$ holds for the bounded component $H$ of
$\Om\setminus[a,b]$.

There is a close  connection between these two notions if
$E\subset\R$.
 \begin{thm}\label{2.2th1} (\cite{andup}) A  set $E\subset \R$ is
uniformly perfect if and only if $\OC\setminus K_E$, defined in
Subsection 2.1, is a John domain.
 \end{thm}
 Since the behavior of a conformal mapping of a John domain onto the
 unit disk is well-studied (see,
 for example, \cite{pom1}), the theorem above can be useful in the
  investigation of metric
properties of the Green function for the complement of a uniformly
perfect subset of $\R$.

 In particular, Theorem \ref{2.2th1} can be used to solve the
 inverse problem of approximation theory of functions that are
 continuous on a uniformly perfect compact subset of the real line (see, for
 details, \cite{andup}).

\subsection{On the Green function for a complement of a compact
subset of $\R$}
 First, we
discuss the following recent remarkable  result by   Totik
\cite{tot}.
 Let $E\subset [0,1]$ be a compact set of positive
logarithmic capacity and let $\Omega$ be the complement of $E$ in
$\CC$. The smoothness of $g_\Om$ and $\mu_E$ at $0$ depends on the
 density of $E$ at $0$. This smoothness can be measured by the function
 $$\theta_E(t):=|[0,t]\setminus E|,\quad t>0.$$
 \begin{thm} (Totik \cite[(2.8) and (2.12)]{tot})\label{tht} There are
 absolute positive constants $C_1,C_2,D_1$ and $D_2$ such that  for $0<r<1$,
 \beq\label{e1.1} g_\Om(-r)\le C_1\sqrt{r}\exp\left(
 D_1\int_r^1\frac{\theta^2_E(t)}{t^3}dt\right)\log\frac{2}{{ \kap}(E)}\,
 ,
 \eeq  \beq\label{e1.2} \mu_E([0,r])\le C_2\sqrt{r}\exp\left(
 D_2\int_r^1\frac{\theta^2_E(t)}{t^3}dt\right).
 \eeq
\end{thm}
 The results in \cite{tot} are formulated and proven for general
compact sets of the unit disk. The theorem above is one of the
main steps in their verification. Even though the statement of
this theorem is rather particular, the theorem  has several
notable  applications, such as Phragm\'en-Lindel\"of-type
theorems, Markov- and Bernstein-type, Remez- and Schur-type
polynomial inequalities, etc.

Observe that we can simplify the geometrical nature of the compact
set $E$ under consideration. Indeed, it is well-known that there
exists a sequence of compact sets $E_n\subset[0,1],\, n\in{\N
}:=\{ 1,2,\cdots\}$, such that

(i) $E\subset E_n$ and each $E_n$ consists of a finite number of
closed intervals,

(ii) for $0<r<1$, we have
$$g_\Om(-r)=\lim_{n\to\infty}g_{\Om_n}(-r),\quad
\Om_n:=\OC\setminus E_n,$$
$$\mu_E([0,r])=\lim_{n\to\infty}\mu_{E_n}([0,r]).$$ The set
$[0,1]\setminus E_n$ is smaller and simpler then $[0,1]\setminus
E$. For example, $$\theta_{E_n}(t)\le\theta_E(t),\quad t>0.$$
However, $g_{\Om_n}$ and $\mu_{E_n}$ can be arbitrarily close to
$g_{\Om}$ and $\mu_{E}$. Thus,  in order to establish Totik-type
results it is natural to concentrate  only on compact sets
consisting of a finite number of real intervals.

Let $$E=\cup_{j=1}^k[a_j,b_j],\quad 0\le a_1<b_1<a_2 <\cdots
<a_k<b_k\le 1,$$ and let $$E^*:=(0,1)\setminus E
=\cup_{j=1}^m(\al_j,\be_j),\quad 0\le \al_1<\be_1<\al_2 <\cdots
<\al_m<\be_m\le 1.$$ For $0<r<1$, we set $E^*_r:=E^*\setminus
(0,r]$. We are interested in the case when $E^*_r\neq\emptyset$,
i.e., $$E^*_r =\cup_{j=1}^{m_r}(\al_{j,r},\be_{j,r}),\quad r\le
\al_{1,r}<\be_{1,r}<\al_{2,r} <\cdots <\al_{m_r,r}<\be_{m_r,r}\le
1.$$
\begin{thm}(\cite{andfni})\label{the1} For $0<r<1$  \beq\label{e2.1} g_\Om(-r)\ge
c_1\sqrt{r}\exp\left(
d_1\sum_{j=1}^{m_r}\frac{\be_{j,r}-\al_{j,r}}{\be_{j,r}}\log
\frac{\be_{j,r}}{\al_{j,r}}\right),\eeq where $c_1=1/16,\,
d_1=10^{-13}$.\end{thm}
 Theorem \ref{the1} provides a lower bound
for the Green function (cf. \cite[(3.5)]{tot}). Since in
(\ref{e2.1}) only the size of the components of $E^*_r$ influences
this bound, one cannot expect to find an upper bound of the same
form. We believe that in a Totik-type theorem not only the size of
the components $(\al_{j,r},\be_{j,r})$ but also their mutual
disposition must be important.

We fix $q>1$. The  set of a finite number of closed intervals
$\{[\de_j,\nu_j]\}_{j=1}^{n} =
\{[\de_j(r,q),\nu_j(r,q)]\}_{j=1}^{n}$, where
$0\le\de_1<\nu_1\le\de_2<\cdots\le\de_n<\nu_n\le 1$, is called a
{\it q-covering} of $E^*_r$ if

(i) $E^*_r\subset\cup_{j=1}^n[\de_j,\nu_j]$,

(ii) either $2\de_j\le\nu_j$, or $q|E^*_r\cap [\de_j,\nu_j]|\le
\nu_j-\de_j.$
\begin{thm}\label{the2} (\cite{andfni}) For $0<r<1,\, q>1$ and any finite q-covering
of $E^*_r$ the inequalities \beq\label{e2.2} g_\Om(-r)\le
c_2\sqrt{r}\exp\left(
d_2\sum_{j=1}^{n}\frac{\nu_{j}-\de_{j}}{\nu_{j}}\log
\frac{\nu_{j}}{\de_{j}}\right)\log\frac{2}{\mb{ cap}(E)}\, ,\eeq
 $$ \mu_E([0,r])\le
c_3\sqrt{r}\exp\left(
d_2\sum_{j=1}^{n}\frac{\nu_{j}-\de_{j}}{\nu_{j}}\log
\frac{\nu_{j}}{\de_{j}}\right)
 $$
  hold with $c_2= 24,\, c_3=5$
and $$ d_2=\max\left(1,\frac{2q^2}{\pi(q-1)^2}\right).$$
\end{thm}
 Notice that the factor $\log(2/$cap$(E))$ on the right of
(\ref{e1.1}) and (\ref{e2.2}) appears only to cover pathological
cases. It is useful to keep in mind that
 $$|E|\le 4 \mb{ cap}(E)\le 1.$$
\begin{cor}\label{core1} (\cite{andfni}) The estimates (\ref{e1.1}) and (\ref{e1.2})
hold with $C_1=384,\, C_2=80$ and $D_1=D_2=120$.\end{cor}
\begin{cor}\label{core2} (\cite{andfni}) For the compact set $$
\tilde{E}:=\{0\}\cup\bigcup_{n=1}^\infty\bigcup_{j=1}^{n^2}\left[
\frac{n^2+j-1}{2^{n+1}n^2},\frac{2n^2+2j-1}{2^{n+2}n^2}\right].$$
we have  $$ g_{\OC\setminus \tilde{E}}(-r)\le c\sqrt{r},\quad
0<r<1,$$ with  some absolute constant $c>0$,  which is better then
(\ref{e1.1}).\end{cor} Indeed, let $$\tilde{E}_r:=\tilde{E}\cap
[r,1],\quad 0<r<1.$$ For $\tilde{E}_r^*=(r,1)\setminus\tilde{E}_r$
with $2^{-k-2}<r\le 2^{-k-1} $ we construct a $2$-covering
$$[r,2^{-k}],\left\{\left\{\left[
\frac{n^2+j-1}{2^{n+1}n^2},\frac{n^2+j}{2^{n+1}n^2}\right]
\right\}_{j=1}^{n^2}\right\}_{n=1}^{k-1},\left[\frac{1}{2},1\right].$$
 By the
monotonicity of the Green function and Theorem \ref{the2}, for any
$0<r<1$ and some absolute constant $c>0$ we obtain $$
g_{\OC\setminus \tilde{E}}(-r)\le g_{\OC\setminus
\tilde{E}_r}(-r)\le c\sqrt{r}.$$

 In what follows in this subsection we assume that $0$ is a
regular point of $E$, i.e., $g_\Om(z)$ extends continuously to $0$
and $g_\Om(0)=0$.

The monotonicity of the Green function yields  $$g_\Om(z)\ge
g_{\OC\setminus [0,1]}(z),\quad z\in \C\setminus [0,1],$$ that is,
if $E$ has the "highest density" at $0$, then $g_\Om$ has the
"highest smoothness" at the origin. In particular,
 \beq\label{2.3ng}
 g_\Om(-r)\ge
g_{\OC\setminus [0,1]}(-r)>\sqrt{r},\quad 0<r<1.
 \eeq
 In
this regard, we would like to explore properties of $E$ whose
Green's function has the ``highest smoothness" at $0$, that is, of
$E$ conforming to the following condition $$g_\Om(z)\le
c|z|^{1/2},\quad c=\mb{const}>0,\,z\in \C,$$ which is known to be
the same as \beq\label{f1.1} \limsup_{r\to 0}
\frac{g_\Om(-r)}{r^{1/2}}<\infty \eeq (cf. \cite[Corollary
III.1.10]{saftot}). Various sufficient conditions for (\ref{f1.1})
in  terms of metric properties of $E$ are stated in \cite{tot},
where the reader can also find further references.

There are compact sets $E\subset[0,1]$ of linear Lebesgue measure
0 with property (\ref{f1.1}) (see e.g. \cite[Corollary 5.2]{tot}),
hence (\ref{f1.1}) may hold, though the set  $E$ is not dense at 0
in terms of linear measure. On the contrary, our first result
states that if $E$ satisfies (\ref{f1.1}) then its density in a
small neighborhood of $0$, measured in terms of logarithmic
capacity, is arbitrarily close to the density of $[0,1]$ in that
neighborhood.
\begin{thm} (\cite{andh})\label{thf11} The condition (\ref{f1.1}) implies
\beq \lim_{r\to 0}\frac{{ \kap} (E\cap[0,r])}{{
\kap}([0,r])}=1.\label{limt}\eeq
\end{thm}
 The converse of Theorem \ref{thf11} is slightly weaker.
\begin{thm} (\cite{andh})\label{thf3} If $E$ satisfies $(\ref{limt})$, then
\beq\label{f1.3}\lim_{r\to 0}\frac{g_{\Om}
(-r)}{r^{1/2-\ve}}=0,\quad 0<\ve<\frac{1}{2}\, .\eeq\end{thm}
 The
connection between properties (\ref{f1.1}), (\ref{limt}) and
(\ref{f1.3}) is quite delicate. For example, even a slight
alteration of (\ref{f1.1}) can lead to the violation of
(\ref{limt}). As an illustration of this phenomenon we construct a
regular set $E\subset[0,1]$ such that $(\ref{f1.3})$ holds  and
\beq\label{f1.4} \liminf_{r\to 0}\frac{\kap(E\cap[0,r])}{
\kap([0,r])}=0.\eeq
 Let
  $$b_j:=2^{-2^{j-1}},\quad
a_j:=b_{j+1}\log(j+1),\quad j\in\N.
 $$
 Consider
 $$E:=\{0\}\cup\left(\cup_{j=1}^\infty[a_j,b_j]\right).$$
  We have
  \beq\label{2.3ng4.1} \lim_{r\to
0}\left(\log\frac{1}{r}\right)^{-1}
\int_r^1\frac{\theta^2_{E}(x)}{x^3}\, dx =0,
 \eeq
 and
  \beq\label{2.3ng4.2}
\lim_{j\to\infty}\frac{b_{j+1}}{a_j}=0. \eeq
 Thus, (\ref{f1.3}) follows from (\ref{e1.1}) and (\ref{2.3ng4.1}).
 Moreover, since
 $$\frac{\kap(E\cap[0,a_j])}{a_j}\le
\frac {b_{j+1}}{4a_j}\, ,$$
 (\ref{2.3ng4.2}) implies (\ref{f1.4}).

A comprehensive description of $E$ satisfying (\ref{f1.1}) was
 recently provided by Carleson and Totik \cite{cartot}.

\subsection{Cantor-type sets}
 Let $0<\ve_j<1$ and $K(j)\in\N,\,
j\in\N$ be two sequences. Starting from $I=[0,1]$ first we remove
$K(1)$ open intervals $I_{1},\cdots,I_{K(1)}$ of $I$ such that
$I\setminus\cup_{k(1)=1}^{K(1)}I_{k(1)}$ consists of $K(1)+1$
disjoint closed intervals $J_{1},\cdots,J_{K(1)+1}$ and
 $$|I_{k(1)}|=\frac{\ve_1}{K(1)},\quad 1\le k(1)\le K(1),$$
  $$|J_{k(1)}|=\frac{1-\ve_1}{K(1)+1},\quad 1\le k(1)\le K(1)+1.$$
 Then, for any $1\le k(1)\le K(1)+1$
 we remove
$K(2)$ open intervals $I_{k(1),1},\cdots,I_{k(1),K(2)}$ of
$J_{k(1)}$ such that
$J_{k(1)}\setminus\cup_{k(2)=1}^{K(2)}I_{k(1),k(2)}$ consists of
$K(2)+1$ disjoint closed intervals
$J_{k(1),1},\cdots,J_{k(1),K(2)+1}$ and
 $$|I_{k(1),k(2)}|=\frac{1-\ve_1}{K(1)+1} \frac{\ve_2}{K(2)}, \quad 1\le k(2)\le K(2),$$
  $$|J_{k(1),k(2)}|=\frac{1-\ve_1}{K(1)+1} \frac{1-\ve_2}{K(2)+1},\quad 1\le k(2)\le K(2)+1,$$
  etc.

  Denote the so obtained Cantor-type set by
  ${\cal{C}}={\cal{C}}(\{\ve_j\},\{K(j)\})$. That is
  ${\cal{C}}:=\cap_{n=1}^\infty
  {\cal{C}}_n,$
  where
  $${\cal{C}}_n={\cal{C}}_n(\{\ve_j\},\{K(j)\}):=\bigcup_{\fk(n)}
  J_{\fk(n)}$$
  is the set we obtain after $n$ steps during the construction,
  and
 $$\fk(j):=k(1),k(2),\cdots,k(j),\quad j\in\N$$
 is a multiple index.
 \begin{thm} \label{2.4th1} (\cite{andcs}) The
 following two conditions are equivalent:

 (i) $g_{\ov{\C}\setminus {\cal{C}}}$ satisfies (\ref{f1.1}) with $E={\cal{C}}$;

 (ii) $\sum_j\ve_j^2<\infty$.
 \end{thm}
 In the case $K(j)=1,\, j\in\N$, this statement is equivalent to
 \cite[Theorem 5.3]{tot}, but the latter is stated for the
 equilibrium measure on ${\cal{C}}$.
Interestingly,  (ii) does not depend
 on the sequence $\{K(j)\}$.

\subsection{On sparse sets with  Green function of the highest
smoothness}
 Let $E\subset\R$ be a compact set  with
positive logarithmic capacity. For simplicity we assume that
$E\subset[-1,1]$ and $\pm1\in E$. Let $\Om=\OC\setminus E$. In
what follows we assume that $E$ is a regular set, i.e., $g_\Om$
extends continuously to $E$ where it takes the value $0$.

We are going to discuss the metric properties of $E$
 such that $g_\Om$ satisfies the $1/2$-H\"older  condition
  \beq\label{2.51.1}
  |g_\Om(z_2)- g_\Om(z_1)|\le c|z_2-z_1|^{1/2},\quad
  z_1,z_2\in\Om\setminus\{\infty\},
  \eeq
 where $c>0$ is some constant.

 According to (\ref{2.3ng})
 the choice of the right-hand side of (\ref{2.51.1}) appears to be
 best suited for this theory.
 In this
regard, we discuss the properties of $E$ whose Green's function
has the ``highest smoothness".

Recently Totik \cite{tot02,tot} constructed two examples of a set
$E$ whose Green's function satisfies (\ref{2.51.1}) and whose
linear measure  is zero.

We analyze the problem: how sparse can  $E$ be, in terms of its
Hausdorff dimension dim$(E)$ \cite[p. 224]{pom1}, if
  it satisfies (\ref{2.51.1}).

First, we note that if $E$ satisfies (\ref{2.51.1}) then
 \beq\label{2.51.2}
 \mb{dim}(E)\ge\frac{1}{2}\, .
 \eeq
 Indeed, from (\ref{2.51.1}) it follows immediately (for details, see
 \cite{cartot}, proof of Proposition 1.4) that for any interval $I\subset\R$,
  $$\mu_E(I\cap E)\le c_1|I|^{1/2},$$
  where $c_1$ is a
  positive constant.

  Hence, for any covering of $E$ by intervals $\{I_j\}\subset\R$
  we have
  $$\sum_j |I_j|^{1/2}\ge c_1^{-1}\sum_j\mu_E(I_j\cap E)\ge c_1^{-1},$$
  which proves (\ref{2.51.2}).
 \begin{thm}\label{2.5th1} (\cite{ands})
 There exists a regular set $E_0\subset\R$ with the following
 properties:

 (i) $g_{\OC\setminus E_0}$ satisfies (\ref{2.51.1});

 (ii) ${\mb{dim}}(E_0)=1/2$.
 \end{thm}
Next, we describe the construction of $E_0$ in Theorem
\ref{2.5th1}.
 For $-1\le a<b\le 1$ we consider two
sequences of real numbers
 $$\cdots<x_{-2}<x_{-1}<x_{0}<x_{1}<x_{2}<\cdots,\quad
 x_k-x_0=x_0-x_{-k}$$
 and
 $$
 y_0>y_{\pm1}>y_{\pm2}>\cdots,\quad y_k=y_{-k},$$
 such that
 $$x_0=\frac{a+b}{2},\quad
 y_0=\frac{b-a}{2}\exp\left\{-\frac{2}{b-a}\right\},$$
 $$y_k=(b-x_k)\exp\left\{-\frac{1}{b-x_k}\right\},\quad
 k\in\N=\{1,2,\cdots\},$$
 $$\frac{y_k}{x_k-x_{k-1}}=\frac{1}{\pi}\left(\frac{1}{b-x_k}-\log
 \frac{1}{b-x_k}\right),\quad k\in\N.$$
 We have
 $$\lim_{k\to\infty}x_{-k}=a,\,\lim_{k\to\infty}x_{k}=b,\,\lim_{k\to\infty}y_{k}=0.$$
 Let $z_k=x_k+iy_k$. For $k\in\Z=\{0,\pm1,\pm2,\cdots\}$ consider vertical intervals
 $J_k=[x_k,z_k]$ and horizontal intervals $I_k=[x_{k-1},x_k]$.
 For multiple indices we use the notation
 $$\fk(m)=k(1),k(2),\cdots,k(m),\quad
 \fk(m)-1=k(1),k(2),\cdots,k(m-1),k(m)-1,$$
  where $m\in\N$ and $k(m)\in\Z$.
  We inductively define  two sets of intervals
  $$\{ J_{\fk(m)}\}_{\fk(m)\in\Z^m}\quad \mb{and}\quad\{
  I_{\fk(m)}\}_{\fk(m)\in{\Z}^m}$$
  in the following way.
 Denote by
 $$\{ J_{\fk(1)}\}_{\fk(1)\in\Z}\quad \mb{and}\quad\{
  I_{\fk(1)}\}_{\fk(1)\in\Z}$$
  the sequences of vertical and horizontal intervals, which we
  obtain by the above procedure  for $[a,b]=[-1,1]$.

  Next, for $m>1$ denote by
  $$\{ J_{\fk(m)}\}_{\fk(m)\in\Z^m}\quad \mb{and}\quad\{
  I_{\fk(m)}\}_{\fk(m)\in\Z^m}$$
   the sequences of vertical and horizontal intervals, which we
  obtain by the above procedure  for $[a,b]=I_{\fk(m-1)}$.
  The endpoints of $\{ J_{\fk(m)}\}$ we denote by $x_{\fk(m)}\in\R$ and
  $z_{\fk(m)}\in\C$ respectively, so that $
  I_{\fk(m)}=[x_{\fk(m)-1},x_{\fk(m)}]$.
  Since
  $$D_0=\{z=x+iy:\,|x|<1,y>0\}\setminus\left(\bigcup_{m\in\N}
  \bigcup_{\fk(m)\in\Z^m}J_{\fk(m)}\right)$$
  is a simply connected domain, by the Riemann mapping theorem
  there exists a conformal mapping $\phi_0$ of $D_0$ onto the
  upper half plane $\He.$

  We interpret the boundary of $D_0$ in terms of Carath\'eodory's
  theory of prime ends (see \cite{pom}).
  Let $P(D_0)$ denote the set of all prime ends of $D_0$.
   For a prime end  $Z\in P(D_0)$ denote its impression by $|Z|$.
  By our construction, all prime ends of $D_0$ are of the first
  kind, i.e., $|Z|$ is a singleton for any $Z\in P(D_0)$.
 For the homeomorphism between $D_0\cup P(D_0)$ and $\ov{\He}$ we
 preserve the same notation $\phi_0$. We
 denote by $\psi_0=\phi_0^{-1}$ the inverse homeomorphism.
 We identify the
 prime end $\psi_0(w), w\in\R$ with its impression when no confusion
  can arise.
   If
  $z\in\partial D_0$ is the impression of only one prime end
 it will also cause no confusion if we use the same letter $z$ to
 designate the  prime end and its impression.
 For example, we write
  $\infty,-1,z_{\fk(m)},1$ for prime ends with impressions
  at those points.

  To define $\phi_0$ uniquely we  normalize it by the boundary
  conditions
  $$\phi_0(\infty)=\infty,\,\phi_0(-1)=-1,\,\phi_0(1)=1.$$
  Each point of $J_{\fk(m)}\setminus\{z_{\fk(m)}\}$ is the
  impression of two prime ends and $z_{\fk(m)}$ is the impression
  of exactly one prime end. Moreover,
  $$\phi_0(\{Z\in P( D_0):\,|Z|\in
  J_{\fk(m)}\setminus\{x_{\fk(m)}\}\})$$
  is an open subinterval of $(-1,1)$ which we denote by
  $J'_{\fk(m)}=(\xi^-_{\fk(m)},\xi^+_{\fk(m)}).$
  Let $\xi_{\fk(m)}=\phi_0(z_{\fk(m)})$.

  In \cite{ands} we show that the compact set
  $$
  E_0=[-1,1]\setminus\left(\bigcup_{m\in\N}
  \bigcup_{\fk(m)\in\Z^m}J'_{\fk(m)}\right)
  $$
  satisfies the conditions of Theorem \ref{2.5th1}. The crucial fact is
  that for $w\in\ov{\He}\cap\Om_0$:
  \beq\label{2.52.2}
  g_{\Om_0}(w)=\frac{\pi}{2}\Im (\psi_0(w)),
  \eeq
  where $\Om_0=\CC\setminus E_0.$

  In order to prove (\ref{2.52.2}), consider the function
  $$h(w)=\left\{\begin{array}{ll}
  \frac{\pi}{2}\Im (\psi_0(w))&\mb{ if }w\in \ov{\He}
  \cap\Om_0,\\[2ex]
  \frac{\pi}{2}\Im (\psi_0(\ov{w}))&\mb{ if }w\in \CC\setminus
  \ov{\He}.
  \end{array}\right.
  $$
  It is continuous in $\Om_0\setminus\{\infty\}$ and, according to
  the distortion properties of $\psi_0$, the difference
  $$h(w)-\log|w|$$
  is bounded in the neighborhood of $\infty$.

  The function $h$ is harmonic in $\C\setminus\R$. In order to
  prove that $h$ coincides with $g_{\Om_0}$ it is sufficient to
  show that $h$ is harmonic in some neighborhood of each
  $$\xi\in(\R\setminus E_0)\setminus \left(\bigcup_{m\in\N}
  \bigcup_{\fk(m)\in\Z^m}\xi_{\fk(m)}\right).
  $$
  Let $\ve=\ve(\xi)>0$ be such that
  $$[\xi-\ve,\xi+\ve]\subset (\R\setminus E_0)\setminus \left(\bigcup_{m\in\N}
  \bigcup_{\fk(m)\in\Z^m}\xi_{\fk(m)}\right).
  $$
  Since all derivatives of $\psi_0$ can be extended continuously
  to $[\xi-\ve,\xi+\ve]$, it is enough to show that for
  $k=1,2;j=0,1,2;j\le k$ and $w=u+iv$:
  $$\lim_{w\to\xi\atop\Im w>0}\frac{\partial^kh(w)}{\partial
  u^j\partial v^{k-j}}=
  \lim_{w\to\xi\atop\Im w<0}\frac{\partial^kh(w)}{\partial
  u^j\partial v^{k-j}}\, ,$$
  which can be easily done.

 It is also natural to consider the problem of how sparse can  set $E$ be such
 that the following local version of (\ref{2.51.1}) is valid:
  \beq\label{2.51.3}
  g_\Om(z)=g_\Om(z)- g_\Om(-1)\le c\,|z+1|^{1/2},\quad
  z\in\Om\setminus\{\infty\},
  \eeq
 where $c>0$ is a constant. The structural properties of compact sets
 satisfying (\ref{2.51.3}) are discussed in \cite{cartot,andh}
 (cf. Subsection 2.3),
 where the density of $E$ near $-1$ is measured in terms of
 logarithmic capacity.
 \begin{thm}\label{2.5th2} (\cite{ands})
  There exists a regular set $E_1\subset\R$ with the following
 properties:

 (i) $g_{\OC\setminus E_1}$ satisfies (\ref{2.51.3});

 (ii) ${\mb{dim}}(E_1)=0$.
 \end{thm}
We describe the construction of $E_1$ in Theorem \ref{2.5th2}.
 We begin with two sequences of real
numbers
 $$1=x_0>x_1>x_2>\cdots>-1\mb{ and } 4=y_0>y_1>y_2>\cdots>0$$
 such that
 $$y_k=(x_k+1)^2,\quad k\in\N,$$
 $$\lim_{k\to\infty}x_k=-1,\quad \lim_{k\to\infty}y_k=0,$$
 $$\frac{y_k}{x_{k-1}-x_k}\ge\frac{2}{\pi}\log
 \frac{1}{x_{k-1}-x_k}\, ,\quad x_{k-1}-x_k<\frac{1}{2}, \quad k\in\N.$$

 Starting with the set of intervals
 $$I_k=[x_{k-1},x_k],\, J_k=[x_{k},x_k+iy_k]=[x_{k},z_k],\,\quad
 k=k(1)\in\N,$$
 we construct the sets of intervals
 $\{I_{\fk(m)}\}$ and $\{J_{\fk(m)}\}$ in the following manner.

 Let for $m\ge 2$, intervals $\{I_{\fk(m-1)}\}$ and $\{J_{\fk(m-1)}\}$
 be constructed, and let
 $$ (A_{\fk(m-1)})^2=
 \exp\left\{m^2+\pi\sum_{j=1}^{m-1}\frac{|J_{\fk(j)}|}{|I_{\fk(j)}|}\right\}.$$
 We define $\de_{\fk(m-1)}>0$ such that
 $$
 \frac{|J_{\fk(m-1)}|}{\de_{\fk(m-1)}}\ge\frac{4m}{\pi}\log
 \frac{A_{\fk(m-1)}}{\de_{\fk(m-1)}}\, .
 $$
 Next, we select a finite number of points
 $$x_{\fk(m-1)-1}=x_{\fk(m-1),0}>x_{\fk(m-1),1}>\cdots>
 x_{\fk(m-1),K(m)}=x_{\fk(m-1)}$$
 such that for any $1\le k(m)\le K(m)$,
 $$
 \frac{1}{2} \de_{\fk(m-1)}\le x_{\fk(m-1),k(m)-1}-x_{\fk(m-1),k(m)}\le
 \de_{\fk(m-1)}.
 $$
 Let
 $$y_{\fk(m)}=\frac{1}{2}y_{\fk(m-1)},\quad z_{\fk(m)}
 =x_{\fk(m)}+iy_{\fk(m)},\quad 0\le k(m)\le K(m),$$
  $$J_{\fk(m)}=[x_{\fk(m)},z_{\fk(m)}],\quad 0\le k(m)\le K(m),$$
 $$I_{\fk(m)}=[x_{\fk(m)},x_{\fk(m)-1}],\quad
 \quad 1\le k(m)\le K(m).$$
 Denote by $\phi_1$ a conformal mapping of the simply connected domain
 $$D_1=\{z=x+iy:\,|x|<1,y>0\}\setminus\left(\bigcup_{m\in\N}
  \bigcup_{1\le k(j)\le K(j)\atop 1\le j\le m}J_{\fk(m)}\right),$$
  where $K(1)=\infty$,
  onto $\He$.

  Let $P(D_1)$ be the set of all prime ends of $D_1$. The reasoning
  about the structure of $P(D_0)$  applies to
  $P(D_1)$.

  We extend $\phi_1$ to
 the homeomorphism  $\phi_1:D_1\cup P(D_1)\to\ov{\He}$ and denote
 the inverse mapping  by $\psi_1=\phi_1^{-1}$.
 Sometimes, for simplicity, we identify  $\psi_1(w), w\in\R$, with
 the impression of $\psi_1(w)$.

 We normalize  $\phi_1$  by the boundary
  conditions
  $$\phi_1(\infty)=\infty,\,\phi_1(-1)=-1,\,\phi_1(1)=1.$$
 For $1\le k(j)\le K(j),1\le j\le m-1$ and $1\le k(m)\le K(m)-1$
   define intervals
  $$J'_{\fk(m)}=(\xi^-_{\fk(m)},\xi^+_{\fk(m)})=
  \phi_1(\{Z\in P( D_1):\,|Z|\in
  J_{\fk(m)}\setminus\{x_{\fk(m)}\}\})$$
  and points $\xi_{\fk(m)}=\phi_1(z_{\fk(m)})$.

  In \cite{ands}, we show that the compact set
  $$
  E_1=[-1,1]\setminus\left(\bigcup_{m\in\N}
  \bigcup_{1\le k(j)\le K(j),1\le j\le m-1\atop 1\le k(m)\le K(m)-1}J'_{\fk(m)}\right)
  $$
  satisfies the conditions of Theorem \ref{2.5th2}. The basic idea is to
  apply the
  formula
  $$
  g_{\Om_1}(w)=\frac{\pi}{2}\Im (\psi_1(w)),\quad
  w\in\ov{\He}\cap\Om_1,
  $$
  where $\Om_1=\CC\setminus E_1$, whose proof is the same as the
  proof of (\ref{2.52.2}).

 We conclude this section with the following remark. One of the
  natural ways to construct sparse sets with H\"older
 continuous Green function is to consider (nowhere dense)
 Cantor-type sets (see \cite{plec,biavol,lit,tot1,tot02}, \cite[Chapter
 5]{tot}).

 Let $\{\ve_j\}$ be a sequence with $0<\ve_j<1$. Starting from
 $[-1,1]$ we first  remove the middle $\ve_1$ part of this
 interval. Then, in the second step, we remove the middle $\ve_2$ part of
 both remaining
 intervals, etc. Denote the so obtained Cantor set by $
{\cal C}= {\cal{C}}(\{\ve_j\})$. According to \cite[Theorem
5.1]{tot} and the
 reasoning in the same monograph \cite[p. 48, after
 Corollary 5.2]{tot}   the following three conditions are equivalent:

 (i) $g_{\ov{\C}\setminus{\cal C}}$ satisfies (\ref{2.51.1});

 (ii) $g_{\ov{\C}\setminus{\cal C}}$ satisfies (\ref{2.51.3});

 (iii) $\sum_j\ve_j^2<\infty$.

 At the same time, by \cite[Theorem 10.5]{pom1} each Cantor type set
 ${ \cal{C}}(\{\ve_j\})$ with the property
 $$\lim_{j\to\infty}\ve_j=0$$
 has  Hausdorff dimension $1$. Therefore, Cantor-type sets
 cannot be used in the proof of either Theorem \ref{2.5th1} or
 Theorem \ref{2.5th2}.

\subsection{Open problems} We begin with
 a new construction of nowhere dense sets.
 It is well-known that Cantor-type sets present  a
 remarkable example of nowhere dense sets which are ``thick" from
 the point of view of potential theory (cf. \cite{car, nev, tot}).
 Motivated by results of this section,
 we suggest the following new construction of such
 sets.
 Let $a_k>0,\, k\in\N$, be such that
 $\lim_{k\to\infty}a_k=0.$
 Starting from the half-strip
 $$\Sigma_0:=\{ z=x+iy:\, |x|<1,y>0\}$$
 we first divide the base $I_0:=[-1,1]$ of $\Sigma_0$ into two
 intervals $I_{1,1}:=[-1,0]$ and $I_{1,2}:=[0,1]$ and  remove
  the vertical slit $J_{1,1}:=[0,ia_1]$  (with one
 endpoint in the middle of  $I_0$). Then, in the second step, we
 divide each of the two new horizontal intervals from the previous step into
 two subintervals of the same length $1/2$ and  remove
  the vertical slits $J_{2,1}:=[-1/2,-1/2+ia_2]$ as well as
   $J_{2,2}:=[1/2,1/2+ia_2]$ (with one
 endpoint in the middle of the base intervals $I_{1,1}$ and
 $I_{1,2}$, respectively), etc.

 As a result, we have a simply connected domain
 $$\Sigma=\Sigma (\{a_k\}):=\Sigma_0\setminus \left(
 \bigcup_{k,m}J_{k,m}\right).$$
 By the Riemann mapping theorem
  there exists a conformal mapping $\phi$ of $\Sigma$ onto the
  upper half plane $\He$.

  We interpret the boundary of $\Sigma$ in terms of Carath\'eodory's
  theory of prime ends (see \cite{pom}).
  Let $P(\Sigma)$ denote the set of all prime ends of $\Sigma$.
  By our construction, all prime ends of $\Sigma$ are of the first
  kind, i.e., $|Z|$ is  singleton for any $Z\in P(\Sigma)$.
 For the homeomorphism between $\Sigma\cup P(\Sigma)$ and $\ov{\He}$,
  which coincides with $\phi$ in $\He$,  we
 preserve the same notation $\phi$.

  To define $\phi$ uniquely we  normalize it by the boundary
  conditions
  $$\phi(\infty)=\infty,\,\phi(-1)=-1,\,\phi(1)=1.$$
  Each interior point of the slit $J_{k,m}=[x_{k,m},x_{k,m}+ia_k]$ is the
  impression of two prime ends. Moreover,
  $$J'_{k,m}:=\phi(\{Z\in P( \Sigma):\,|Z|\in
  J_{k,m}\setminus\{x_{k,m}\}\})$$
   is an open subinterval of $(-1,1)$.

   Hence,
   $$E=E(\{a_k\}):=[-1,1]\setminus\left(\bigcup_{k,m}J'_{k,m}\right)$$ is
   a
   nowhere dense subset of $[-1,1]$.

    It seems to be an interesting problem  to investigate the connection between
    the  geometry of
    $E$ (for example, its Hausdorff dimension and  Hausdorff
    measure), the rate of decrease of $a_k$ as $k\to\infty$, and
    continuous properties of the Green function $g_{\OC\setminus
    E}$.

    The crucial fact is
  that for $w\in\ov{\He}\cap\Om$:
  $$
  g_{\Om}(w)=\frac{\pi}{2}\Im (\phi^{-1}(w)),
  $$
  where $\Om=\OC\setminus E.$

  For example,  the following problems
  can be considered.

   {\bf Problem 1.} {\it Are the
 following two conditions

 (i) $g_{\Om}$ satisfies the  the  $1/2$-H\"older
 property, i.e.,
 $$g_\Om(z)\le
 c\mb{  dist}(z,E)^{1/2},\quad z\in \Om,$$
 where $c=c(E)>0$ is a constant and
 $$\mb{ dist}(A,B):=\inf_{\z\in A,\z\in B}|z-\z|,\quad
 A,B\subset{\C};$$

 (ii) $\sum_ja_j^2<\infty$

  equivalent?}

  (cf.  \cite[Theorem
5.1]{tot}  concerning Cantor-type sets).

  {\bf Problem 2.} {\it Use the ideas of this section to streamline the
  proof of the  Carleson - Totik  \cite[Theorem 1.1]{cartot}
  characterization of compact sets $E\subset\R$ such that the Green function
  $g_{\OC\setminus E}$
  satisfies the H\"older condition, i.e., there are constants
  $c>0$ and $0<\alpha\le1/2$ such that}
 $$g_{\OC\setminus E}(z)\le
 c\mb{ dist}(z,E)^{\alpha},\quad z\in \C\setminus E.$$
  We conjecture that a more general choice of
 horizontal intervals $I_{k,m}$ and  slits $J_{k,m}$ in the
 procedure described above will allow one to construct
 nowhere dense sets with various extremal properties.

 Consider a typical example. Let $h(r),0\le r\le1/2$, be a
 monotone increasing function and $h(0)=0$. Denote by
 $\Lambda_h(E)$ the Hausdorff measure of a set $E\subset\C$ with
 respect to $h$ (see \cite[p. 224]{pom1}). A well-known metric
 criterion for sets of zero capacity  states that (see \cite[Theorem
 3.14]{lan}) if
 $$\Lambda_h(E)<\infty,\quad h(r)=|\ln r|^{-1},$$
 then cap$(E)=0$.

 {\bf Problem 3.} {\it Show that for any  monotone
 increasing function $g(r),0\le r\le1/2$,
 satisfying
  $$\lim_{r\to 0}\frac{g(r)}{h(r)}=0,$$
 there exists a compact set $E_g\subset\R$ such that}
 $$\mb{cap}(E_g)>0\quad\mb{ and}\quad\Lambda_g(E_g)<\infty.$$
 (cf. \cite[Chapter IV]{car}).

\section{Remez-type inequalities}

\subsection{Remez-type inequalities in terms of capacity}
 Let $\bP_n$ be the set of all real polynomials of degree at most
 $n\in\N.$
 The
Remez inequality \cite{rem} (see also \cite{erd,borerd,gan})
asserts that
 \beq \label{c1.1}
 ||p_n||_I\le
T_n\left(\frac{2+s}{2-s}\right)
 \eeq
 for every
$p_n\in\bP_n$  such that
 \beq \label{c1.11}
 |\{
x\in I:\, |p_n(x)|\le 1\}|\ge 2-s,\quad 0<s<2,
 \eeq
 where
$I:=[-1,1]$,  $T_n$ is the Chebyshev polynomial of degree $n$, and
$||\cdot||_A$ means the uniform norm along $A\subset\C$.

Since $$ T_n(x)\le (x+\sqrt{x^2-1})^n,\quad x>1, $$ we have by
(\ref{c1.1}) that a polynomial $p_n$ with (\ref{c1.11}) satisfies
\beq \label{c1.12}
||p_n||_I\le\left(\frac{\sqrt{2}+\sqrt{s}}{\sqrt{2}-\sqrt{s}}\right)^n.
\eeq
 The last inequality (more precisely
its $n$-th root) is asymptotically sharp.

Remez-type inequalities give bounds for classes of functions on a
line segment, on a curve or on a region of the complex plane,
given that the modulus of the functions is bounded by $1$ on some
subset of prescribed measure. Remez-type inequalities play a
central role in proving other important inequalities for
generalized nonnegative polynomials, exponentials of logarithmic
potentials and M\"untz polynomials.
 There are a
number of recent significant advances in this area.
 A survey of results concerning various
generalizations and numerous applications of this classical
inequality can be found in \cite{erd}, \cite{borerd} and
\cite{gan}. In particular,
 a pointwise, asymptotic version of
(\ref{c1.1}) is also obtained  \cite[Theorem 4]{erd92}. Namely
 \beq\label{3.1p4}
|p_n(x)|\leq
\exp\left(c\,n\,\min\left\{\frac{s}{\sqrt{1-x^2}},\sqrt{s}\right\}\right)
 \eeq
  holds for $x\in I$ and every $p_n\in \bP_n$  satisfying (\ref{c1.11}),
  where $c>0$ is some
universal constant.

 In this section, we discuss an analogue of (\ref{c1.11})--(\ref{c1.12})
in which we use logarithmic capacity instead of linear length. Our
main results deal not only with polynomials, but also with
exponentials of potentials (see \cite{erd, erdlisaf}).

Given a nonnegative Borel measure $\nu$ with compact support in
$\C$ and finite total mass $\nu(\C)>0$ as well as a constant
$c\in\R$, we say that $$ Q_{\nu,c}(z):=\mb{exp}(c-U^\nu(z)),\quad
z\in\C, $$ where $$ U^\nu(z):=\int\log\frac{1}{|\z-z|}\,
d\nu(\z),\quad z\in\C ,$$ is the logarithmic potential of $\nu$,
is an {\it exponential of a potential} of degree $\nu(\C)$.

Let $$ E_{\nu,c}:=\{ z\in \C:\, Q_{\nu,c}(z)\le 1\}. $$ Theorem
2.1 and Corollary 2.11 in \cite{erdlisaf} assert that for $0<s<2$
the condition \beq \label{c1.2} |E_{\nu,c}\cap I|\ge 2-s \eeq
implies \beq \label{c1.3}
||Q_{\nu,c}||_I\le\left(\frac{\sqrt{2}+\sqrt{s}}{\sqrt{2}-\sqrt{s}}
\right)^{\nu(\C)}.
 \eeq
 \begin{thm}(\cite{andrl})\label{thc1} Let $0<\de<1/2$. Then the condition \beq \label{c1.4}
\mb{ cap} (E_{\nu,c}\cap I)\ge \frac{1}{2}-\de \eeq yields that
\beq \label{c1.5}
||Q_{\nu,c}||_I\le\left(\frac{1+\sqrt{2\de}}{1-\sqrt{2\de}}
\right)^{\nu(\C)}. \eeq \end{thm}
 Since $|E_{\nu,c}\cap I|\le 4\,\kap (E_{\nu,c}\cap I)$ \cite[p.
337]{pom}, the assertion (\ref{c1.2})--(\ref{c1.3}) follows from
(\ref{c1.4})--(\ref{c1.5}).

Furthermore, for $0<\de<1/2$, set $$
\nu=\nu_\de:=\mu_{[-1,1-4\de]},\quad
c=c_\de:=\log\frac{2}{1-2\de}. $$
 Then $\nu(\C)=1$, $
E_{\nu,c}=[-1,1-4\de], $
 $$ Q_{\nu,c}(x)=\frac{1}{1-2\de}\left(
x+2\de+((x+2\de)^2-(1-2\de)^2)^{1/2} \right),\quad x\ge 1-4\delta.
$$ Therefore, in this case $$ \kap (E_{\nu,c}\cap
I)=\frac{1}{2}-\de, $$ $$
||Q_{\nu,c}||_I=Q_{\nu,c}(1)=\frac{1+\sqrt{2\de}}{1-\sqrt{2\de}},
$$ which shows the sharpness of Theorem \ref{thc1}.

Note that the modulus of any complex polynomial
$p_n(z)=c\prod_{j=1}^{n}(z-z_j),\,0\neq c\in\C$ can be written as
an exponential of a potential in the following way. Let
 \beq\label{3.1p1}
\nu_n:=\sum_{j=1}^n\de_{z_j},
 \eeq
  where $\de_z$ is the Dirac unit measure in a point $z\in\C$.
For $z\in\C$, we have
 \beq\label{3.1p2}
 Q_{\nu_n,\log |c|}(z)=\mb{exp }(\log
|c|+\log \prod_{j=1}^{n}|z-z_j|) =|p_n(z)|.
 \eeq
  Therefore, applying
 Theorem \ref{thc1} we obtain for $0<\delta<1/2$: the condition $$ \kap
(\{ x\in I:\, |p_n(x)|\le 1\})\ge \frac{1}{2}-\de $$ implies $$
||p_n||_I\le\left( \frac{1+\sqrt{2\de}}{1-\sqrt{2\de}} \right)^n
$$ (cf. (\ref{c1.11})--(\ref{c1.12})).

 The previous remark can be rewritten in a form as
in \cite[Theorem 1.1]{cuydrilub}. Namely, let $r>0$ and
$p_n\in\bP_n$ be  such that $||p_n||_{[-r,r]}=1$. Then for
$0<\ve<1$, $$ \kap (\{ x\in [-r,r]:\, |p_n(x)|\le
\ve^n\})\le\frac{2r\ve}{(1+\ve)^2}. $$ This inequality is
asymptotically sharp for any fixed $\ve$ and $r$.

Next, we present an analogue of the above results  for  complex
polynomials. By $\Pn$ we denote the set of all complex polynomials
of degree at most $n\in\N$. Let $$ \Pi (p_n):=\{ z\in
\C:|p_n(z)|>1\} ,\quad p_n\in\Pn . $$

From the numerous  generalizations of the  Remez inequality, we
cite one result which is a direct consequence of the trigonometric
version of the Remez inequality (and is equivalent to this
trigonometric version, up to constants).

Assume that $p_n\in\Pn$, $\T:=\{ z:|z|=1\}$ and \beq \label{d1.1}
| \T\cap\Pi (p_n)|\leq s,\quad 0<s\le\frac{\pi}{2} . \eeq Then,
$q_n(t):=|p_n(e^{it})|^2$ is a trigonometric polynomial of degree
at most $n$ and, by the Remez-type inequality on the size of
trigonometric polynomials (cf. \cite[Theorem 2]{erd92}, \cite[p.
230]{borerd}), we have
 \beq \label{d1.2} ||p_n||_{\T}\leq e^{2sn},
\quad 0<s\leq\frac{\pi}{2}.
 \eeq
 Our next
objective is to discuss an analogue of (\ref{d1.1})--(\ref{d1.2})
in which we use logarithmic capacity instead of linear  length. As
before, our main result deals not only with polynomials, but also
with exponentials of potentials.
\begin{thm} (\cite{andrc}) \label{thd1} Let $0<\de<1 $. Then the condition $$
 {\kap}(E_{\nu,c}\cap \T) \ge  \de $$  implies that $$
 ||Q_{\nu, c}||_{\T} \le \left(\frac{1+\sqrt{1-\de^2}}{
\de } \right)^{\nu(\C)}. $$\end{thm}
 In order to examine the sharpness of Theorem
\ref{thd1} we consider the following example.

Let $0<\al<\pi/2$, and let
 \beq\label{3.1p3}
 L=L_\al:=\{e^{i\theta}:\,
2\al\le\theta\le 2\pi-2\al\}. \eeq
 Since the function
$$z=\Psi(w)=-w\frac{w-a}{1-aw},$$ where $a=1/\cos\al$, maps
$\Delta=\CC\setminus\OD$ onto $\Omega:=\OC\setminus L$ (cf.
\cite{goo}) and since the Green function of $\Om$ with pole at
$\infty$ can be defined via the inverse function $\Phi:=\Psi^{-1}$
by the formula $$g_\Om(z)=\log |\Phi(z)|,\quad z\in\Omega,$$ we
have \beq\label{d1.15} \kap
(L)=\lim_{w\to\infty}\frac{\Psi(w)}{w}=\frac{1}{a}=\cos\al,\eeq as
well as
\begin{eqnarray} \max_{z\in\T\setminus L}g_\Om(z)&=&
g_\Om(1) =\log |\Phi(1)|\nonumber\\ \label{d1.16}
&=&\log(a+\sqrt{a^2-1})= \log\frac{1+\sqrt{1-\kap( L)^2}}{\kap
(L)}.\end{eqnarray}

Let $c=c_\al:=-\log \kap( L)$ and let $\nu=\nu_\al:=\mu_L$ be the
equilibrium measure for $L$; that is, $\nu(\C)=1$. Since for
$z\in\C$, $$ U^\nu(z)=-g_{\OC\setminus L}(z)-\log\kap( L),$$ and
therefore $$Q_{\nu,c}(z)=\mb{exp }(g_{\OC\setminus L}(z)),$$ we
have $E_{\nu,c}=L$ as well as
$$||Q_{\nu,c}||_{\T}=\frac{1+\sqrt{1-\kap( L)^2}}{\kap (L)}.$$
This shows the exactness of Theorem \ref{thd1}.

Applying Theorem \ref{thd1} to the exponential of a potential
defined by (\ref{3.1p1})--(\ref{3.1p2}), we obtain the following:
 for $p_n\in\Pn$ the condition \beq \label{d1.11}
\kap( \T\setminus\Pi(p_n)) \ge \de,\quad 0<\delta<1
 \eeq
  yields
\beq\label{d1.12} ||p_n||_{\T}\le
 \left(\frac{1+\sqrt{1-\de^2}}{
\de } \right)^{n}. \eeq
 Since for any $E\subset\T$ we have $\kap
(E)\ge\sin\frac{|E|}{4}$ (see \cite{pom1}),
(\ref{d1.11})--(\ref{d1.12}) imply the following refinement of
(\ref{d1.1})--(\ref{d1.2}): For $p_n\in\Pn$ the condition
\beq\label{d1.21} | \T\cap\Pi(p_n)| \le s,\quad 0<s<2\pi, \eeq
implies \beq\label{d1.22} ||p_n||_{\T}\le
 \left(\frac{1+\sin\frac{s}{4}}{\cos\frac{s}{4}}
  \right)^{n}.\eeq This result is also sharp in the following
  sense. Let $0<s<2\pi$, $\alpha=s/4$, and let $L=L_\al$ be
  defined as in (\ref{3.1p3}). By (\ref{d1.15}) and
  (\ref{d1.16}), $$
  g_\Om(1)=\log\frac{1+\sin\frac{s}{4}}{\cos\frac{s}{4}}\,
  .$$
 We denote by $f_n(z)$ the $n$-th Fekete
  polynomial for a compact set $L$ (see \cite{saftot}). Hence,  condition
  (\ref{d1.21}) holds for  the polynomial
  $p_n(z):=f_n(z)/||f_n||_L$. At the same time,
since
$$\lim_{n\to\infty}\left(\frac{|f_n(z)|}{||f_n||_L}\right)^{1/n}=
  \exp (g_\Om(z))
  ,\quad z\in\Om\setminus\{\infty\}$$
  (see \cite[p. 151]{saftot}), we have
  $$\lim_{n\to\infty}|p_n(1)|^{1/n}=\frac{1+\sin\frac{s}{4}}{\cos\frac{s}{4}}\,
  $$
  (cf. (\ref{d1.22})).

\subsection{Remez-type inequalities in the complex plane}
 Let
$ m_2(A)$ be   the two-dimensional Lebesgue measure (area) of a
set $A\subset\C$. The analogue of (\ref{c1.1}), where the unit
interval $[-1,1]$ is replaced by the closure $\ov{G}$ of some
bounded Jordan domain $G\subset \C$ and (\ref{c1.11}) by
\begin{equation}
\label{3.2eq1.3}
m_2\left(\left\{z\in\overline{G}\,:|p_n(z)|\leq1\right\}\right)\geq
m_2(\overline{G})-s,\quad 0<s<m_2(\overline{G}),
\end{equation}
is studied by Erd\'elyi, Li, and Saff \cite{erdlisaf}. Let
$\Pn(\overline{G},s)$ denote the subset of polynomials in $\Pn$
satisfying (\ref{3.2eq1.3}), and let
 $$
R_n(z,s):=\sup_{p_n\in\Pn(\overline{G},s)}|p_n(z)|,\quad z\in
L:=\partial G.
 $$
  If $L$ is a $C^2$-curve it is established in
\cite{erdlisaf} that there is a constant $c_j=c_j(G)>0$ where $
j=1,2$, such that
 \beq \label{3.2eq1.4} R_n(z,s) \leq\mb{exp}(c_1n\sqrt{s}),
\quad z\in L,\quad 0<s\leq c_2<m_2(\overline{G}).
 \eeq
 Actually, this result is established in a
more general context of exponentials of logarithmic potentials,
where it is used to prove Nikolskii-type inequalities (cf.
\cite{erd},\cite{erdlisaf}). The same problem  was investigated
recently \cite[Theorem 2.3]{kro}  for domains with smooth boundary
(under weaker restrictions on the smoothness rate than in
\cite{erdlisaf}).

We generalize the above  results in two directions: we  obtain
pointwise bounds for $R_n(z,s)$, depending on $z\in L$, and we
replace the strong $C^2$ restriction for the boundaries of $G$ by
weaker ones. Our results can easily be generalized to exponentials
of logarithmic potentials as well. The method to obtain our (sharp
up to constants) estimates differs from the approaches used
elsewhere \cite{erdlisaf},\cite{kro}. We make use of properties of
Green's functions (cf. \cite{saftot}), principles of
symmetrization (cf. \cite{be}), and the technique of moduli of
families of curves (cf. \cite{ah}, \cite{lv}), combined with a
useful estimate from \cite{ar}.

To aid in further discussion we introduce additional notation. For
$z\in\C$ and $ r>0$ let
  $$
D(z,r):=\left\{\zeta \,:\,|z-\zeta |<r\right\}, \quad
  D(r):=D(0,r),$$
 $$
C(z,r):=\left\{\zeta \,:\,|z-\zeta |=r\right\}, \quad
  C(r):=C(0,r).$$
 Let $G\subset
\C$ be a bounded Jordan domain, $\overline{\C} :=\C \cup
\{\infty\}$ and
 $$ L:=\partial G,\quad \Omega
:=\overline{\C}\setminus\overline{G}, \quad \gamma_z(r):=\Omega
\cap C(z,r). $$

 We use the convention that $c_1,c_2,\dots$ denote positive
constants and $\ve_1,\ve _2,\dots$ sufficiently small positive
constants. If  not stated otherwise, we assume that both types of
constants depend only on $G$.
 \begin{thm} (\cite{arcp})
\label{3.2th1} Let $G$ be a bounded domain and $z\in L:=\partial
G$. Suppose that there are constants $\ve _1,\ve _2$ such that
\begin{equation}
\label{3.2eq1.5} \ve _1r\le |\gamma _{z}(r)|\le (2\pi -\ve
_1)r,\quad 0<r<\ve _2.
\end{equation}
 Then there exist constants $c_1,c_2,c_3>c_2$ and $\ve
_3<\left(\frac{\ve_2}{2c_3}\right)^2$ depending only on $\ve _1$
and $ \ve _2$ such that
 \begin{equation}
\label{3.2eq1.6}
  R_n(z,s)
 \leq \exp\left(c_1n \exp\left(-\pi
\int\limits_{c_3\sqrt{s}}^{\ve _2/2}
\frac{dr}{h_{z,c_2\sqrt{s}}(r)}\right)\right),  \quad 0<s<\ve _3,
\end{equation}
where
  $$h_{z,\de}(r):=\sup_{|t-r|\leq\de}|\ga_z(t)|,\quad 0<\de<r.$$
 Moreover, for any arbitrary bounded domain $G$
\begin{equation}
\label{3.2eq1.61}
R_n(z,s)\leq\left(\frac{c_4}{\MG-s}\right)^{n},\quad 0< s<\MG,
\end{equation}
 where $c_4>2\,\MG$ depends only on the diameter of $G$.
\end{thm}
 The inequality (\ref{3.2eq1.6}) is the main statement of the theorem. The
 estimate (\ref{3.2eq1.61}) is included for the completeness of the
 result.
 The condition (\ref{3.2eq1.5}) excludes any domain with a cusp at
 $z$.

 In the proof of Theorem \ref{3.2th1} we exploit the following deep connection between
estimates which express the possible growth of a polynomial with a
known norm on a given compact set $E\subset\C$ and the behavior of
the Green function for $\Om=\ov{\C}\setminus E$.

For $z\in \Omega $ and $u>0$, the following two conditions are
equivalent:

(i) $ g_\Om(z)\le u$;

(ii) for any $p\in\Pn$ and $n\in\N$, $$ |p(z)|\leq e^{un}||p||_E.
$$

Indeed, (i) $\Rightarrow$ (ii) follows from the Bernstein--Walsh
lemma \cite[p. 77]{wal}. (ii) $\Rightarrow$ (i) is a  consequence
of a result by Myrberg and Leja (see  \cite[p. 333]{pom}).

We  study the properties of the Green function by methods of
geometric function theory (using symmetrization, moduli of curve
families, distortion theorems, harmonic measure, etc.) which allow
us, according to the implication (i) $\Rightarrow$ (ii), to get
(\ref{3.2eq1.6}).

Note that the sharpness of the results  for the Green function
means, by virtue of the equivalence (i) and (ii), the sharpness
(up to  constants) of the corresponding Remez-type inequalities.

Since for an arbitrary Jordan domain $G$, $z\in L$ and $0<\de<r$,
 $$
  h_{z,\de}(r)\le 2\pi(r+\de),
  $$
 according to (\ref{3.2eq1.6}), for each domain satisfying
(\ref{3.2eq1.5}) we obtain
 \beq \label{3.2eq1.11}
  R_n(z,s)\le
\mb{exp}(c_{5}ns^{1/4}),\qquad 0<s<\ve_3 .
 \eeq
 An example below shows that for domains with cusps at $z$ the
inequality (\ref{3.2eq1.11})  in general does not hold (i.e., the
restriction (\ref{3.2eq1.5}) in Theorem \ref{3.2th1} is
essential).

Indeed, let $k>2$ be fixed and let
 \beq\label{3.2aa.7}
  G=G_k:=\left\{
z=x+iy:0<x<1,0<y<\frac{1}{2}x^{k-1}\right\} .
 \eeq
 For $1/2<\delta
<1$ we have
 $$ m_2(\overline G\setminus D(1,\delta ))\le
\frac{1}{2}\int\limits_0^{2(1-\delta )}x^{k-1}
dx=\frac{2^{k-1}}{k}(1-\delta )^k. $$
 For
$0<s<1/(2k)=m_2(\overline{G_k})$ let
 $$ \delta =\delta
(s):=1-\left( \frac{ks}{2^{k-1}}\right)^{1/k}. $$ The polynomial
$$ p_n(z):=\left( \frac{z-1}{\delta}\right)^n $$ belongs to
${\Pa}_n(\ov{G},s)$. However, $$ |p_n(0)|=\frac{1}{\delta^n}\ge
\mb{exp }(\ve_4ns^{1/k}). $$
 Hence,
for $k>4$ and  $G$ defined by (\ref{3.2aa.7}), the inequality
(\ref{3.2eq1.11}) is violated.

If more information is known about the geometry of the domain $G$,
the expression in (\ref{3.2eq1.6}) can be made more explicit. The
following example illustrates this point. A Jordan curve  is
called {\it Dini-smooth } (cf. \cite[p. 48]{pom1}) if it is smooth
and if the angle $\beta (s)$ of the tangent, considered in terms
of the arclength $s$, satisfies $$ \left|\beta (s_2)-\beta
(s_1)\right|<h(s_2-s_1),\quad s_1<s_2, $$ where $h(x)$ is an
increasing function for which
 \beq
\label{3.2eq1.63} \int\limits_{0}^{1}\frac{h(x)}{x}dx<\infty .
\eeq We call a Jordan arc Dini-smooth if it is a subarc of some
Dini-smooth curve.
 \begin{thm} (\cite{arcp})
\label{3.2th2} Let $L=\partial G$ consist of finitely many
Dini-smooth arcs $l_j$, which form exterior angles $\alpha _j\pi
,\ 0<\alpha _j<2,$ at their junction points $z_j, j=1,\dots,m$.
Let $z\in L$ be arbitrary and let $z_k$ be the nearest point to
$z$ among the $z_j$, i.e., $$ |z_k-z|=\min_{1\leq j\leq m}|z-z_j|.
$$ Then, for $0<s<\MG$ the inequality
\begin{equation}
\label{3.2eq1.7} R_n(z,s)\leq \exp\left(
\frac{c_6\,n\sqrt{s}}{(\sqrt{s}+|z_k-z|)^{1-1/\alpha
_k}}\log\left(\frac{c_4}{\MG-s}\right)\right)
\end{equation}
 holds.
\end{thm}
 Note that Theorem \ref{3.2th2} implies (and extends)
(\ref{3.2eq1.4}) since  we can choose $\alpha _1=\al_2=1$ and
arbitrary distinct points $z_1,z_2\in L$.

Now, let us proceed with the discussion of the sharpness of
 (\ref{3.2eq1.7}).

Let $\Phi$ denote the Riemann function that maps $\Omega$
conformally and univalently onto $\Delta :=\overline{\C}\setminus
\overline{\D}$, where $\D:=\{ z:\,|z|<1\}$ is the unit disk, and
is normalized by the conditions $$ \Phi (\infty )=\infty ,\quad
\Phi '(\infty )>0. $$
 We extend $ \Phi$ to
the homeomorphism $\Phi :\overline{\Omega}\rightarrow
\overline{\Delta}$. For $z\in {\C} $ and $\delta >0$ let $$
L_\delta :=\{ z\in \Omega :|\Phi (z)|=1+\delta\} , $$
 $$ \rho_\delta
(z):=\mb{dist}(z,L_\delta ). $$
 Let the function $\delta (z,t)$ be defined by the relation
 $$
\rho_{\delta (z,t)}(z)=t,\quad z\in L,t>0. $$
 Observe  that under the
assumptions of Theorem \ref{3.2th2}
 $$ \delta\asymp\frac{\rho_\delta
(z)}{ (\rho_\delta (z)+|z-z_k|)^{1-1/\alpha_k}}, \quad 0<\delta
<1, $$ where $a\asymp b$ denotes the double inequality
$\ve_5\,b\le a\le c_7\,b$.

Indeed, let $w:=\Phi (z),w_k:=\Phi (z_k),w_\delta :=(1+\delta )w,
z_\delta :=\Phi^{-1}(w_\delta )$. According to the distortion
properties of conformal mappings of domains with piecewise
Dini-smooth boundary (cf. \cite[p. 52]{pom1} or \cite[p.
33]{andbla})  we have $$ \rho_\delta (z)\asymp |z-z_\delta |\asymp
\delta (\delta +|w-w_k|)^{\alpha_k -1}\asymp \delta (\rho_\delta
(z)+|z-z_k|)^{1-1/\alpha_k}. $$
 Hence (\ref{3.2eq1.7}) is equivalent to
  \beq \label{3.2eq1.8}
 R_n(z,s)\leq \exp\left(
c_8\, n\,
\de(z,\sqrt{s})\log\left(\frac{c_4}{\MG-s}\right)\right), \,
0<s<\MG. \eeq
 In \cite{arcp} the sharpness of (\ref{3.2eq1.8}) is established for an
arbitrary quasidisk $G$, i.e., a Jordan domain bounded by a
quasiconformal curve $L:=\partial G$.

Recall that by  Ahlfors' theorem (see, for example, \cite[p.
100]{lv}) a Jordan curve $L$ is {\it quasiconformal}  if and only
if there exists a constant $c_9$ such that for $z_1,z_2\in L$
 \beq\label{3.3ah}
 \mb{min }\{ \mb{diam}(L'), \mb{diam}(L'')\} \le c_9|z_1-z_2|,
 \eeq
  where $L'$ and
$L''$ denote the two components of $L\setminus\{ z_1,z_2\} $.
 Thus, we exclude from our consideration the regions with cusps on
the boundary.

Using Ahlfors' criterion, one can easily verify that convex
curves, curves of bounded variation without cusps, and rectifiable
Jordan curves which have the same order of arc length and chord
length are quasiconformal. At the same time, each part of a
quasiconformal curve can be nonrectifiable.

  The domain from Theorem \ref{3.2th2} is a quasidisk.

There is a simple way to show that there exists $\ve_6$ such that
for
 $$
m_2(\overline{G})-\ve_6<s<m_2(\overline{G}) $$
 the inequality (\ref{3.2eq1.8}) is  sharp (up to constants).
 Indeed, let $L=\partial G$ be quasiconformal, $z\in L$. Denote by
$z^*\in L$ an arbitrary point satisfying
 $$ |z-z^*|\ge \frac{1}{2}
\mb{ diam }(\ov{G}). $$
 Choose $0<r<|z-z^*|$ such that $$
m_2(D(z^*,r)\cap G)=\MG -s $$ (this is always possible if $ \ve_6$
is sufficiently small). Since $G$ is a quasidisk, for any $z\in L$
 (\ref{3.2eq1.5}) holds . Therefore,
  $$ r\le c_{10}(\MG -s)^{1/2} $$
  holds as well.

  The polynomial $$ p_n(\z
):=\frac{(\z -z^*)^n}{r^n} $$
 belongs to $\Pa_n(\ov{G},s)$. At
the same time
 $$
 R_n(z,s)\ge |p_n(z)|\ge \left(\frac{\ve_7}{\MG
-s}\right)^{n/2}=\exp\left(\frac{n}{2}\log\frac{\ve_7}{\MG
-s}\right),
 $$
 which shows the sharpness of (\ref{3.2eq1.8})
 for values of $s$ close to $\MG$.

 If $0<s\le \MG-\ve_6$, then (\ref{3.2eq1.8}) implies
 $$\limsup\limits_{n\to\infty}\, \frac{\log R_n(z,s)}{n\delta (z,\sqrt{s}) }
 \leq c_{11}=c_{11}(\ve_6),\quad z\in L. $$
 \begin{thm} (\cite{arcp})
\label{3.2th3} Let $G$ be a quasidisk. There exists a constant
$\ve_8=\ve_8(\ve_6)$
  such that for $z\in L$ and $0<s\le \MG-\ve_6$ the inequality
 \beq\label{3.2aa.3}
 \liminf\limits_{n\to\infty}\, \frac{\log R_n(z,s)}{n\delta (z,\sqrt{s})
 }\ge\ve_8
 \eeq
 holds.
\end{thm}
 The inequality (\ref{3.2aa.3}) demonstrates that (\ref{3.2eq1.8}) is
 asymptotically sharp (with respect to $n$)  for all $z\in L$ and
 $0<s\le \MG-\ve_6$.

 \subsection{Pointwise Remez-type inequalities in the unit disk}
 Observe that (\ref{3.2eq1.4}) states a uniform bound.
  Our  objective is to derive the pointwise extension of
  this bound.
  Note that the pointwise extension (\ref{3.1p4}) of the
  classical Remez inequality (\ref{c1.1}) was established relatively recently
  by Erd\'elyi \cite[Theorem 4]{erd92}. The proof of this theorem is based on a
  Remez-type inequality for trigonometric polynomials (cf. \cite[p. 230]{borerd}, \cite{gan}).
  Our approach is quite different.
   We use ideas from
  potential theory in the plane \cite{ran}, \cite{saftot}, principles of
  circular symmetrization
  \cite{be},
  and estimation of conformal invariants such as moduli of families of curves
   \cite{ah}, \cite{lv}.

For simplicity, we formulate the appropriate result only for the
unit disk.
 \begin{thm}\label{3.3th1} (\cite{andud})
 The condition
 $$ m_2({\ov{\D}}\cap E_{\nu,c})\ge\pi-s,\quad 0<s<\pi,
 $$
  yields
 $$ Q_{\nu,c}(z)\le\exp\left(c_1\nu(\C)\sqrt{s}
 \exp\left(-c_2\frac{(1-|z|)^2}{s}\right)\right),\quad z\in\OD,
 $$
 where  $c_1$ and $c_2$ are  positive absolute constants.
 \end{thm}
 Theorem \ref{3.3th1} is sharp in the following
sense.
 Given $0<s<\pi/2$ and $0\le x<1$. Let $\de:=1-x$. Set for $0<r<1$,
 \beq\label{3.1p5}
 \tilde{E}_{r,\de}:=\OD\setminus\left(\{\z:\,|\z-x|<r\}
 \cup\{\xi+i\eta:\,x<\xi\le 1,|\eta|<r\}\right).
 \eeq
 Since
 $$m_2(\tilde{E}_{r,\de})>\pi-\frac{\pi r^2}{2}-2\de r,$$
 taking $r$ such that $(\pi r^2)/2+2\de r=s$, i.e.,  $$
 r=r(\de,s):=\frac{s}{\de+\sqrt{\de^2+(\pi/2)s}}\, ,$$
 we have for $\tilde{E}:=\tilde{E}_{r,\de}$,
 \beq\label{3.31.5} m_2(\tilde{E})>\pi-s.\eeq
 In \cite{andud} it is shown that
 \beq\label{3.31.6}
 g_{\OC\setminus\tilde{E}}(x,\infty)\ge
 c_3\sqrt{s}\exp\left(-4\frac{(1-x)^2}{s}\right),\quad
 c_3=\frac{e^{-10\pi}}{2\sqrt{2}}.\eeq

 Let $$\nu=m\,\mu_{\tilde{E}},\quad c=-m\,\log\kap(\tilde{E}),$$
 where $m>0$ is an arbitrary number.
 Since
 $$U^{\mu_{\tilde{E}}}(z)=\left\{\begin{array}{ll}
 -  g_{\OC\setminus\tilde{E}}(z,\infty)-\log\kap(\tilde{E}), &
 z\in\C\setminus\tilde{E},\\[2ex]
 -\log\kap(\tilde{E}), &
 z\in\tilde{E},\end{array}\right.$$
 we have
 $$ \nu(\C)=m,\quad E_{\nu,c}=\tilde{E}.$$
 Moreover, (\ref{3.31.6}) implies
 \beq\label{3.31.8}
 Q_{\nu,c}(x)\ge\exp\left(c_3\,m\,\sqrt{s}
 \exp\left(-4\frac{(1-x)^2}{s}\right)\right).\eeq
 Relations (\ref{3.31.5}) and (\ref{3.31.8}) show the
 sharpness of Theorem \ref{3.3th1} for $z\in\D$ (its sharpness
 for $z\in\T$ is known \cite{erdlisaf}).

 Applying Theorem
\ref{3.3th1} to the exponential of a potential defined by
(\ref{3.1p1})--(\ref{3.1p2}), we obtain the following: for any
complex polynomial $p_n\in\Pn$ the condition
 $$
 m_2(\{z\in\OD:\,
|p_n(z)|\le 1\})\ge \pi-s,\quad 0<s<\pi $$
 implies
 $$ |p_n(z)|\le \exp\left(c_1\, n\,\sqrt{s}
 \exp\left(-c_2\frac{(1-|z|)^2}{s}\right) \right),\quad z\in\OD.
 $$
 The last inequality  is  sharp for large $n$.
 That is, given $0<s<\pi/2$ and $0\le x<1$, let
 $\tilde{E}$ be defined as in (\ref{3.1p5}).
 It was shown in \cite{andud} that for
 $$ n>\frac{8\log 2}{c_3\sqrt{s}}
 \exp\left(4\frac{(1-x)^2}{s}\right)
 $$
 there exists a polynomial $P_n\in\Pn$ such that
 $$|P_n(z)|\le 1,\quad z\in\tilde{E},
 $$
 $$ |P_n(x)|\ge \exp\left(\frac{c_3}{2}\, n\,\sqrt{s}
 \exp\left(-4\frac{(1-x)^2}{s}\right)\right).
 $$

\subsection{Remez-type inequalities in terms of linear measure}
Next, we discuss analogues of (\ref{c1.1}) and (\ref{d1.2}) with
an arbitrary Jordan arc or curve instead of $[-1,1]$, and  a
quasismooth (in the sense of Lavrentiev)
 curve instead of $\T$,
respectively.
 Our  results
deal not only with polynomials but also with exponentials of
logarithmic potentials (cf. \cite{erd, erdlisaf}).

 Let $L\in\C$ be a bounded Jordan arc or curve.
 For a (Borel) set
$V\subset L$ we consider its covering $U=\cup_{j}U_j\supset V$ by
a finite number or countably many  open (i.e., without endpoints)
subarcs $U_j$ of $L$.
 Let
  $$
\sigma_L(V):=\inf\sum\limits_{j}\mb{diam}(U_j),
 $$
  where the
infimum is taken over all such open coverings of $V$.

Note that
 \beq\label{3.4aa.3}
 \sigma_L(V)\le \min\{|V|,\mb{ diam}(L)\}.
 \eeq
 For an exponential of a potential $Q_{\nu,c}$ set
 $$
 E^*_{\nu,c}:=\C\setminus E_{\nu,c}=\{z\in\C:\, Q_{\nu,c}(z)>1\}.
 $$
 \begin{thm} (\cite{arl})
\label{3.4th2} Let $L$ be an arbitrary  Jordan arc or curve,
 and let
 $$ \frac{\sigma_L(E_{\nu,c}^*\cap L)}{\mb{ diam}(L)}=:u<\frac{1}{2}.
 $$
  Then
 \beq \label{3.4eq2.6}
||Q_{\nu,c}||_L\le\left(\frac{1+\sqrt{2u}}{1-\sqrt{2u}}\right)^{\nu(\C)}.
\eeq
\end{thm}
 Theorem \ref{3.4th2} extends \cite[Theorem 2.1]{erdlisaf} from the case
 where
 $L=[-1,1]$ to the case where $L$ is an arbitrary Jordan arc or curve.

 Theorem \ref{3.4th4} and the left-hand side of (\ref{3.4eq2.10a}) below show that
 (\ref{3.4eq2.6})
is sharp
 (with respect to the degree $1/2$ of $u$) even for the case of
Jordan curves. However, if we take into consideration additional
information about the
 geometry of $L$, the estimate (\ref{3.4eq2.6})
can be improved.

Let $L$ be a
 {\it quasismooth}
(in the sense of Lavrentiev) curve which is defined by the
following condition. For any $z_1,z_2\in L$,
 \beq\label{3.4aa.1}
  \mb{min }\{
 |L'|, |L''|\}\le c_1|z_1-z_2|,\quad c_1=c_1(L)\ge 1,
 \eeq
where $L'$ and $L''$ are the connected components of $L\setminus\{
z_1,z_2\}$.

According to (\ref{3.4aa.3}) and (\ref{3.4aa.1}), for any
quasismooth curve $L$ and a (Borel) set $V\subset L$ we have
 $$
 \ve_1 |V|\le\sigma_L(V)\le |V|.
 $$
 We proceed with the case where $E$ is a {\it Lavrentiev domain}, i.e.,
$L=\partial E$ is quasismooth.

Denote by $\Phi$  the  conformal mapping of $\Omega=\CC\setminus
E$ onto the exterior $\Delta := \overline{\C} \backslash
\overline{\D}$ of the unit disk $\D$  normalized by the conditions
   $$
   \Phi(\infty) = \infty,\quad  \Phi'(\infty) := \lim_{z \to \infty}
   \frac{\Phi(z)}{z} > 0.
   $$
 Let
 $$
  L_\delta :=\{ z\in\Omega :|\Phi(z)|=1+\delta\} ,\quad
\delta >0,
 $$
  and let the function $\delta (t,L),\, t>0,$  be defined by the relation
 $$ \mb{dist}(L,L_{\delta (t,L)})=t.
 $$
 \begin{thm} (\cite{arl})
\label{3.4th3} Let $L$ be a quasismooth curve and suppose that
 $$
  |E^*_{\nu,c}\cap L|\le s<\frac{1}{2}\, \mb{ diam}(L).
  $$
 Then
  $$
||Q_{\nu,c}||_L\le \exp (c_2\,\delta(s,L)\, \nu(\C))
 $$
 holds with $c_2=c_2(L).$
\end{thm}
In order to discuss the sharpness of the bound  of Theorem
\ref{3.4th3} we consider an important particular case of
exponentials of potentials. Let $V\subset L$ consist of a finite
number of open subarcs of $L$ whose closures are disjoint,
 $ J:=L\setminus V$, and let $c=c(V):=-\log \kap(J)$.  Since
 $$
U^{\mu_J}(z)=-g_{\CC\setminus J}(z)-\log\kap (J),\quad z\in\C,
 $$
 we have
 $$
 Q_{\mu_J,c}(z)=\mb{exp}(g_{\CC\setminus J}(z)), \quad
E^*_{\mu_J,c}\cap L=V. $$
 For $0<s<\mb{diam}(L)$ we set
$$ {{\cal{U}}(s,L)}:=\{ V\subset L: \, V=\cup_{j=1}^mV_j,V_j\mb{
is an open arc },\ov{V_j}\cap\ov{V_k}=\emptyset,
\sum_{j=1}^m\mb{diam}(V_j)\le s\} ,
 $$
  $$
   \la
(s,L):=\sup\limits_{V\in{\cal{U}}(s,L)} \sup\limits_{z\in
V}g_{\CC\setminus J}(z).
 $$
 For any quasismooth curve $L$,  Theorem \ref{3.4th3} implies
   that
  \beq \label{3.4eq2.81} \la (s,L)\le c_2\,\delta (s,L),\quad
   0<s<
\frac{1}{2}\mb{ diam}(L) .
   \eeq
 \begin{thm}
\label{3.4th4} Let $L$ be a quasismooth curve.
  Then
   \beq \label{3.4eq2.82}
   \la (s,L)\ge\ve_2 \delta
(s,L) ,\quad  0<s< \mb{ diam}(L)
 \eeq
 holds with $\ve_2=\ve_2(L)$.
 \end{thm}
 Theorem \ref{3.4th3} and inequalities (\ref{3.4eq2.81}) and  (\ref{3.4eq2.82}) show
 that  the growth of the
 exponentials of logarithmic potentials
  can be related to the function $\delta (s,L)$ which depends only on
  the geometry of $L$. Below we state
some remarks concerning its properties.

By the Ahlfors criterion (\ref{3.3ah}), any quasismooth curve is
quasiconformal. Therefore, $\Phi$ can be extended to a
quasiconformal homeomorphism $\Phi:\CC\to\CC$.
 Taking into account  the
distortion properties of conformal mappings with a quasiconformal
extension (cf. \cite[pp. 289, 347]{pom}), for any quasismooth
curve $L$ we have
 $$ \delta (s,L)\le c_3s^\alpha
,\quad 0<s<\mb{diam}(L),
 $$
  where $c_3=c_3(L)$ and $\alpha
=\alpha (L)>1/2$.
 Thus, for any particular quasismooth curve $L$,
 Theorem \ref{3.4th3} presents better estimates
  (with respect to the order of
 $s$)  than Theorem \ref{3.4th2}.

Next, we introduce the notion of  Dini-convex curves.
 In the remainder of this section we assume  that $L$  is a
 quasismooth
curve. The set $\ov{\C}\setminus L$ consists of two Jordan
domains:
 a bounded one
 $G:=$ int$(L)$ and an unbounded one $\Omega :=$ ext$(L)$.
Let $h$ be a positive nondecreasing function satisfying the
Dini-condition (\ref{3.2eq1.63}), and let for $0<\ve\le 1$
   $$
   W(h,\varepsilon) := \{ \, \z = re^{i\theta}:\, 0 < r < \varepsilon,\,
   \pi h(r) < \theta
   < \pi (1 - h(r))\, \} .
   $$
 We say that  $ L$ is {\it Dini-convex} with respect to $G$ if there
exist $0 < \varepsilon=\ve(L) \le 1$ and a function $h=h_L$
satisfying (\ref{3.2eq1.63}) such that $h(\ve)< 1/2$ and for any
$z \in L$ $$
   \{\z = z+e^{i\theta} \xi: \xi \in W(h,\varepsilon)\} \subset G
$$
 holds with some $0 \le \theta = \theta(z) < 2\pi$.

For example, if there is  $0<r<1$ depending only on $L$ such that
for each $z \in L$ there exists an open disk $D_z$ with radius $r$
such that $D_z \subset  G$ and $\overline{D}_z \cap L = \{z\}$,
then $L$ is Dini-convex with respect to $G$ (with $h(x) = c_4\,x$
and $ \ve=r$).
 \begin{thm} (\cite{arl})
\label{3.4lem1} Let $L$ be a quasismooth curve which is
Dini-convex with respect to G. Then
 \beq \label{3.4eq2.10} \ve_3s\le \delta
(s,L)\le c_5s,\quad 0<s<\mb{ diam}(L).
 \eeq
\end{thm}
 Comparing Theorem \ref{3.4th3} with the right-hand side of
(\ref{3.4eq2.10})  we obtain an analogue of (\ref{d1.2}) for
curves $L$ instead of the unit circle $\T$ (see also
\cite{andej}). This result is sharp because of Theorem
\ref{3.4th4} and the left-hand side of (\ref{3.4eq2.10}).

If $L$ consists of a finite number of Dini-smooth arcs which meet
in the angles $\alpha_1\pi ,\cdots ,\alpha_m\pi$ (with respect to
$\Omega$), $0<\alpha_j<2$, $$ \alpha :=\max (1,\alpha_1 ,\cdots
,\alpha_m), $$ then according to the distortion properties of a
conformal mapping of a domain with piecewise Dini-smooth boundary
onto a  disc (cf. \cite[Chapter 3]{pom1}) we have
 \beq \label{3.4eq2.10a}
\ve_4s^{1/\alpha}\le \delta (s,L)\le c_6s^{1/\alpha},\quad
0<s<\mb{diam}(L).
 \eeq
 Notice that Theorem \ref{3.4th3} and the right-hand side of (\ref{3.4eq2.10a})
 imply  a new Remez-type inequality for domains
 with a piecewise Dini-smooth boundary which is sharp because of
 Theorem \ref{3.4th4} and the left-hand side of (\ref{3.4eq2.10a}).

 \subsection{Open problems}
 We conjecture that the
following sharp pointwise extension of (\ref{c1.1}) and
(\ref{3.1p4}) is valid.

Let $0<s<2$ and $-1<x<1$. We define $a=a(s,x)$ and $b=b(s,x)$ such
that:

(i) $-1<a<x<b<1$ and $b-a=s$;

(ii) the conformal mapping $F=F_E$ defined by (\ref{fdef}) for
$E=[0,\frac{a}{2}+\frac{1}{2}]\cup[\frac{b}{2}+\frac{1}{2},1]$
maps $\frac{x}{2}+\frac{1}{2}$ into the point
 $w_0=\max_{w\in K_E}|w|$.

 Let $\tilde{T}_{n,s}(\xi)=\xi^n+\cdots\in\bP_n$ be the {\it Chebyshev-Akhiezer
 polynomial} deviating least from zero on $E$, that is,
 $$||\tilde{T}_{n,s}||_E=\min_{p\in\bP_{n-1}}||\cdot^n+p(\cdot)||_E.$$
 Let
 $T_{n,s}(\xi):=\tilde{T}_{n,s}(\xi)/||\tilde{T}_{n,s}||_E$.

 {\bf Problem 4.} {\it Is it true that the inequality
 $$|p_n(x)|\le |T_{n,s}(x)|$$
 holds for every $p_n\in\bP_n$ satisfying (\ref{c1.11})?}

 Our next problem concerns the  Remez-type inequality for
 polynomials on a quasidisk.

 {\bf Problem 5.} {\it Let $G$ be a quasidisk, i.e., a Jordan
 domain bounded by a quasiconformal curve $L$.
 Is it true that for $z\in L$ and arbitrary sufficiently small
 positive
 constant $\ve$
  the inequalities
   $$\exp(c_1 n\delta (z,\sqrt{s}))\le R_n(z,s)\le
  \exp(c_2  n\delta (z,\sqrt{s}))
 $$
  hold for any  $0<s\le
m_2(\overline{G})-\ve$ with some constants}
$c_j=c_j(G,\ve),j=1,2?$

\section{Polynomial Approximation}

\subsection{Approximation on an unbounded interval} We consider
functions $f$, continuous and real valued on the non-negative real
line $\R ^+$ and possessing also the
 properties
 \beq\label{b1.1} f>0 \:
\mb{ on } \: \R^+,\quad \lim _{x\to\infty} f(x)=\infty \: . \eeq

For every positive integer $n \in \N$, we define
 $$
 \rho_n(f):=\inf_{p_n\in \bP_n} || \frac{1}{f} -
 \frac{1}{p_n} ||_{\R^+} \: .
 $$
In the present section, we  discuss  necessary and sufficient
conditions for the geometric convergence of reciprocals of
polynomials to the reciprocal of the function $f$ on
 $\R^+$,
  i.e., the inequality
 \beq \label{b1.3}\limsup_{n\to \infty }\rho_n(f)^{1/n}
=\frac{1}{q}< 1 \:. \eeq

The first results in this area were due to  Cody,  Meinardus and
 Varga \cite{codmeivar} concerning the function $ \exp (x)$.
Later,  Meinardus and  Varga \cite{meivar} extended these results
to the class of entire functions of completely regular growth. The
paper \cite{meiredtayvar} gave rise to investigations devoted to
enlarging the class of functions that admit geometric
approximation by reciprocals of polynomials on $\R^+$.

We introduce some notations. Given two numbers $r>0$ and $s>1$,
denote by ${\cal E}_r(s)$ the closed ellipse with foci at the
points $x = 0$ and $x = r$ such that the ratio between the
semimajor axis and semiminor axis equals $(s^2 + 1)/(s^2 - 1)$.

The following theorem states
 remarkable  necessary conditions for
geometric convergence.
 \begin{thm}\label{thbA} (Meinardus \cite{mei}, Meinardus, Reddy,
Taylor, Varga \cite{meiredtayvar}) Let $f$ satisfy (\ref{b1.3}).
Then

 (i) the function $f$ can be extended from $\R^+$
 to an entire function of finite order;

(ii) for every  number $s>1$, there exist  positive constants $c
_1=c_1(s,q)$, $\theta =\theta (s,q)$  and $r_0=r_0(s,q)$ such that
the inequality \beq\label{b1.4} || f||_{{\cal E} _r(s)} \le
c_1||f||_{[0,r]}^{\theta} \eeq holds for all $r\ge r_0$.
\end{thm}
 After the appearance of \cite{meiredtayvar}, a lot of work was
done to find sufficient conditions for (\ref{b1.3}) (cf.
\cite{bla2}--\cite{blakov},  \cite{redshi}, \cite{routay},
\cite{henrou}). The most general known result in this direction is
the following statement.
\begin{thm}\label{thbB} (Blatt, Kovacheva \cite{blakov})
  Assume that $f$ is an entire function
with (\ref{b1.1}) and,
 in addition to condition (\ref{b1.4}),
the inequality \beq \label{b1.5}|| f ||_{[0,r]} \le \mu
(r)^{\lambda},
 \eeq
  where
   $\mu
(r):= \min _{x\ge r}\{f(x)\}$,
  holds for some number  $\lambda > 1$ and for
every $r>r_0$. Then (\ref{b1.3}) is true. \end{thm}
 On  the other
hand,  Henry and  Roulier  \cite{henrou} have shown that the
conditions (i) and (ii) of Theorem \ref{thbA} are not sufficient
for geometric convergence. For example, in \cite{henrou} it was
proved that \beq\label{b1.6} f(x)=1+x+e^x\sin ^2x \eeq
 cannot be approximated with geometric speed.
Their proof was based on the fact that $f$ satisfying (\ref{b1.3})
cannot oscillate too often.

The main goal of this section is to discuss a new necessary and a
new sufficient condition for geometric convergence found in
\cite{andblakov}.

We begin with a necessary condition. Let $f$ be as above, i.e.,
$f$ is an entire function with (\ref{b1.1}). For $ r>0$, we define
the set $$ Z_r:=\{\,  0<x<\infty :f(x)<r \, \}  \: . $$ Then $Z_r$
is the union of a finite number of disjoint open intervals. This
follows from (\ref{b1.1}) and the uniqueness theorem for analytic
functions. Now, we consider the closure $\bar Z_r$ of $Z_r$, which
is regular and possesses a Green's function
 $$g_r(z):=g_{\CC\setminus \bar{Z_r}}(z)$$
  with respect
to the region $\overline{\C} \setminus \bar Z_r$ with pole at
infinity, where  $ g_r := 0 $ on  $ \bar Z_r$.
 For $s>1$ we denote by
${\cal E}_r(f,s)$ the set which consists of the interior of the
level set of $g_r(z)$ and the level set itself for a fixed
parameter $s$, i.e., $$ {\cal E}_r(f,s):= \{ \, z\in \C \, : \,
0\le g_r(z)\le \log s \, \} \:. $$ Then the new necessary
condition for geometric convergence can be formulated as follows.
 \begin{thm}(\cite{andblakov})\label{thb1}
Let $f$ satisfy (\ref{b1.3}). Then for every $1<s<q$ there exist
positive constants $c=c(s,q)$, $\theta =\theta (s,q)$ and
$r_0=r_0(s,q)$ such that \beq\label{b2.1} ||f||_{{\cal
E}_r(f,s)}\le c \, r^\theta, \quad r\ge r_0 \:.\eeq
\end{thm}
 Next, we are going to discuss the geometrical meaning of condition (\ref{b2.1}).
 For $ \infty > H>h>\min _{x\in \R^+}f(x)>0$,
 we introduce the strip domain
 $$
 S(h,H):=\{ \, (x,y) \, : \, -\infty <x<\infty , \: h<y<H \,\}
 $$
 as well as the intersection of this strip with the graph of $f$,
 i.e.,
 $$
 Y(f,h,H):=S(h,H) \cap \{ \, (x,y) \, : \, x\ge 0, \: y=f(x) \,\} \:,
 $$
and define $N(f,h,H)$  to be the number of connected components of
$Y(f,h,H)$ joining the line $\{\Im (z)=h\} $ with the line $\{\Im
(z)=H\} $. Since $f$ satisfies (\ref{b1.1}),  the number
$N(f,h,H)$ is finite and, moreover, it is odd.
\begin{thm}(\cite{andblakov})\label{thb2}
Let $f$ be an entire function satisfying (\ref{b1.1}). If, in
addition, for some $s>1$ and $\theta >1$, the function  $f$
satisfies (\ref{b2.1}), then, for each $M>\theta $,
\beq\label{b2.2} \limsup_{h\to \infty }\frac{N(f,h,h^M)}{\log
h}<\infty \:. \eeq
\end{thm}
 Note that the result of Theorem \ref{thb2} is sharp in the following sense:
For each $M>1$ there exists an entire function $f=f_M$ which
satisfies (\ref{b2.1}) with some $s>1$ and $1<\theta <M$ and
\beq\label{b2.3} \limsup_{h\to \infty }\frac{N(f,h,h^M)}{\log h}>0
\:. \eeq

Indeed, consider the function $$ f_M(x):= e^x+e^{2Mx}\sin ^2\pi x
\:. $$ Obviously, it satisfies the conditions of Theorem
\ref{thbB}. Therefore, $f_M$ guarantees the geometrical
convergence of best approximants in the sense of (\ref{b1.3}),
and, by Theorem \ref{thb1}, $f$ satisfies (\ref{b2.1}) in which we
can take $s$ so close to 1 that $\theta <M$.
 The relation (\ref{b2.3}) immediately follows, if we set $h=e^k, \: k\in \N ,$
and let $k\to \infty $.

The new sufficient condition for geometrical convergence of best
approximants can be stated in the following form.
\begin{thm}(\cite{andblakov})\label{thb3}
Let $f$ be an entire function satisfying (\ref{b1.1}) and
(\ref{b2.1}) with some $s>1$ and $\theta >1$. In addition, assume
that there exists a constant $M=M(f)>1$ such that
 $$ \limsup_{h\to
\infty }N(f,h,h^M)<\infty \:. $$ Then $f$ satisfies (\ref{b1.3}).
\end{thm}
 It is easy to see that Theorem \ref{thbB} follows from Theorem
\ref{thb3}, because under the assumptions of Theorem \ref{thbB},
for $M>\lambda$ and $h$ sufficiently large, we have
$N(f,h,h^M)=1$. At the same time, condition (\ref{b1.4}) is weaker
than (\ref{b2.1}). The  example of the function (\ref{b1.6})
 shows that conditions (\ref{b1.4}) and (\ref{b2.1}) are not equivalent.
Indeed, $f$ given by (\ref{b1.6}) obviously satisfies
(\ref{b1.4}). On the other hand, some
 straightforward calculations show that relation (\ref{b2.2}) does not hold for this function.
Thus, $ f $ does not possess (\ref{b2.1}).

The fact that Theorem \ref{thb3} is essentially stronger than
Theorem \ref{thbB} is not so obvious.
\begin{thm}(\cite{andblakov})\label{thb4}
There exists an entire function $f$ satisfying the assumptions of
Theorem \ref{thb3}, but not possessing property (\ref{b1.5}).
\end{thm}
 The proof of Theorem \ref{thb3} is  based on an analogue of the
classical result due to Bernstein, concerning polynomial
approximation of functions analytic in the neighborhood of a
subinterval of the real axis, for the case of several intervals.

Let $E = \bigcup^k_{j=1} I_j$  be the union of $k$ disjoint
intervals
 $I_j = [\alpha_j,
\beta_j] $ of the real axis $ \R$ and let $\Omega := \overline{\C}
\backslash E$.  The set
   $$
   E^s := \{z \, \in \Omega \, : \, g_\Omega(z) = \log s \,\}, \quad s >
   1,
   $$
consists of at most $k$ (mutually exterior) curves. Denote by $
\mbox{ext} (E^s)$ the unbounded component of $\overline{\C}
\backslash E^s$ and set $ \mb{int}(E^s) := \overline{\C}
\backslash \overline{\mb{ext} (E^s)}$.
 Denote by $C(E)$ the class of all real functions continuous on
 $E$.
 \begin{thm}(\cite{andblakov})\label{thb5}
For each  $f\in C(E)$ satisfying the following two conditions:
 \beq\label{b2.5}
 \mb{ for some } s>1, \, \: f \,\mb {can be extended analytically into }
 \overline{\mb{int}(E^s)} \:,
\eeq \beq\label{b2.6}
 f \: \,\mb{ has at least one zero on each } \: I_j \:,
\eeq
 there exist constants $q >1$ and $c
> 0$ depending only on $s$ and $k$ such that
    \beq\label{b2.7}
    \inf_{p_n\in\bP_n}||f-p_n||_E=:E_n(f,E) \le c \, ||f||_{E^s}
    \, q^{-n}, \quad n \in \N \:.
    \eeq
\end{thm}
 Note that (\ref{b2.7}) can be interpreted as a result concerning
geometric convergence of the polynomials of best approximation to
the function $f$, independent of the geometry of $E$.

The proof of Theorem \ref{thb5} is based on a new concept of
 Faber-type polynomials
for $E$ which can be described as follows.
 Let $E$ be as defined in Theorem \ref{thb5}.
   Denote by
$ g_\Om(z,z_0), \: z,z_0 \in \Omega :=\overline{\C}\setminus E$,
the Green function for $\Omega$ with
 pole at $z_0$. It has a multiple--valued harmonic conjugate $\tilde{g}_\Om(z,z_0)$.
 Thus, the analytic function
    $$
    \Phi(z,z_0)  :=  \exp(g_\Om(z,z_0) + i \tilde{g}_\Om(z,z_0))
    $$
is also multiple-valued.

  Let $\Phi(z):=\Phi(z,\infty)$ and let $n \in \N$ be  arbitrary. If $\Phi(z)^n$ is  single--valued
in $\Omega$, we set
   $$
   W_n(z) := \Phi(z)^n, \quad z \in \Omega \:.
   $$
If $\Phi(z)^n$ is not single-valued in $\Omega$, then according to
\cite[pp. 159, 227]{wid},
 there exist $q \le k-1$ points $x_{1,n},\cdots,x_{q,n} \in
[\alpha_1,\beta_k] \backslash E$ such that the function
   $$
   W_n (z) := \frac{\Phi(z)^n}{\Pi^q_{i=1} \Phi(z,x_{i,n})}, \quad z \in \Omega
   $$
 is single-valued in $\Omega$.
 In both cases the function
    $$
    F_n(z) := \frac{1}{2\pi i} \; \int_{\Gamma} \; \frac{W_n(\zeta)}
    {\zeta-z} d\zeta, \quad z \in \mb{int}(\Gamma) ,
    $$
    where $\Gamma\in\C\setminus E$ is any curve surrounding $E$,
 is a polynomial of degree $ n$, which coincides with the classical Faber
 polynomial in the case of connected $E$.

\subsection{The Nikol'skii-Timan-Dzjadyk-type theorem}
 Let
$E\subset\R$ be a compact set, and let $\om(\de),\de>0$, be a
function of modulus of continuity type, i.e., a positive
nondecreasing function with $\om(0+)=0$ such that for some
constant $c\ge 1$, $$ \om(t\de)\le c\, t\,\om(\de),\quad
\de>0,t>1. $$ Let $C_\om(E)$ consist of all $f\in C(E)$ such that
$$ |f(x_1)-f(x_2)|\le c_1\,\om(|x_2-x_1|),\quad x_1,x_2\in E, $$
with some  $c_1=c_1(f)>0$.

For $\om(\de)=\de^\al,0<\al\le 1$, we set $C_\om(E)=:C^\al(E)$.

One of the central problems in approximation theory is to describe
the relation between the smoothness of functions and the rate of
decrease of their approximation by polynomials when the degree of
these polynomials tends to infinity. The following well-known
statement is the starting point of our consideration.
 \begin{thm} (Nikol'skii \cite{nik}, Timan \cite{tim}, Dzjadyk
\cite{dzj56})\label{thntd} Let $f\in C([-1,1])$ and let $\om$ be a
function of modulus of continuity type satisfying the inequality
\beq \label{a1.1} \de\,\int\limits_\de^1\frac{\om(t)}{t^2}\, dt\le
c_2\, \om(\de),\quad 0<\de<1, \eeq with some constant $c_2>0$.
Then the following assertions are equivalent:

(i) $f\in C_\om([-1,1])$;

(ii) for any $n\in\N$ there exists $p_n\in\bP _n$ such that the
inequality \beq \label{a1.2} |f(x)-p_n(x)|\le
c_3\,\om\left(\frac{1}{n^2}+\frac{\sqrt{1-x^2}}{n} \right),\quad
-1\le x\le 1, \eeq holds with some constant $c_3>0$. \end{thm}
 In
the late 50s - early 60s Dzyadyk \cite{dzj59}, \cite{dzj62} laid
the foundation for a new constructive theory of functions on
continua in the complex plane (a survey of the results and a
bibliography can be found in the monographs \cite{tam},
 \cite{dzj},  \cite{gai},
\cite{she}, \cite{andbeldzj}). He used the following simple but
fundamental idea.

Denote by $I_{1/n}, n\in\N$, the ellipse with foci at $\pm 1$ and
sum of semiaxes equal to $1+1/n$. Such an ellipse is the image of
the circle $\{ w:\, |w|=1+1/n\}$ under the conformal mapping
$z=\frac{1}{2}(w+\frac{1}{w})$ of $\De:=\{ w:\, |w|>1\}$ onto
$\OC\setminus [-1,1]$,  i.e., $I_{1/n}$ is the level line of the
conformal mapping $$ \Phi(z)=z+\sqrt{z^2-1} $$ of $\OC\setminus
[-1,1]$ onto $\De$, where the square root is chosen so that
$\Phi(z)=2z+O(\frac{1}{|z|})$ in a neighborhood of $\infty$.

Then for $-1\le x\le 1$ and $n\in\N$, $$
\frac{1}{n^2}+\frac{\sqrt{1-x^2}}{n}\asymp\rho_{1/n}(x), $$ where
 $$ \rho_{1/n}(x):=
\mb{dist}(x,I_{1/n}). $$
 The concepts of $C_\om,\Phi,I_{1/n}$ and
$\rho_{1/n}(x)$ are
 also meaningful for an
arbitrary bounded continuum in the complex plane. This is the key
to a generalization of the Nikol'skii-Timan-Dzjadyk theorem to
classes of functions on continua in $\C$.

 If $E\subset\C$ is a compact set, then the
interpretation of the Nikol'skii-Timan-Dzjadyk theorem above
 can be rephrased by consideration of the
Green function $g_\Omega$  and its level lines.
 The case when $\Om =\OC\setminus E$ is multiply
connected is discussed  in  \cite{ngu, ngu88, shi1, shi99, mezshi,
andntd}. Each time the extension of a result from the case of a
continuum to the case of a compact set uses quite specific and
non-trivial constructions.

In \cite{and05} we  found how, in the case of finitely connected
$\Om$, the extension of the Nikol'skii-Timan-Dzjadyk theorem can
be obtained
 by using the well-known Bernstein-Walsh lemma on the growth of
a polynomial outside the compact set and  Walsh's  theorem on
polynomial approximation of a function analytic in a neighborhood
of a compact set with connected complement. Our approach is based
on the following theorem.
  \begin{thm} (\cite{and05})
 \label{4.33th2}  Let $E=\cup_{j=1}^mE_j$ consist of $m\in\N,\, m\ge 2,$ disjoint
continua $E_j,\, f\in A(E),\, ||f||_E\le 1$, and let
$z_1,\cdots,z_N\in E$ be distinct points. Let for any $n>
n_0\in\N$ and $j=1,\cdots,m$ there be a polynomial $p_{n,j}\in\Pn$
such that
 $$|f_j(z)-p_{n,j}(z)|\le \ve_j\left(\frac{1}{n},z\right),\quad z\in \partial E_j,$$
 $$
p_{n,j}(z_l)=f_j(z_l),\quad z_l\in E_j ,$$
 where $f_j:=f|_{E_j}$ is the restriction of $f$ to $E_j$,
 and the function $\ve_j(\de,z),\, 0<\de\le 1,\, z\in \partial
 E_j$, satisfies, for any $j=1,\cdots,m$ and $z\in\partial E_j$, the
 properties:

 (i) $\ve_j(\de,z)$ is monotonically increasing in $\de$;

 (ii) $|\ve_j(\de,z)|\le 1,\quad \de\le\de_0\le 1.$

 Then for any $n\in \N,\, n> c_4(n_0+1/\de_0)$ there exists a
 polynomial $p_n\in\Pn$ such that
$$ |f(z)-p_n(z)|\le \ve_j\left(\frac{c_5}{n},z\right)+c_6\,
e^{-c_7n},\quad z\in
\partial E_j,\, j=1,\cdots,m, $$
$$ p_n(z_l)=f(z_l),\quad l=1,\cdots,N,
 $$
  where $c_k,\,
k=4,5,6,7$, depend only on $E$ and the choice of points
$z_1,\cdots,z_N$.
 \end{thm}
 The case of infinitely connected $\Om$ is extremely difficult to
handle. This can be seen from a recent paper by Shirokov
\cite{shi99}.

 In what follows we are going to discuss the case $E\subset\R$,
where the number of components of $E$ can be infinite. It turns
out that the appropriate analogue of the Nikol'skii-Timan-Dzjadyk
theorem is valid for some $E$ that are not "too scarce" (see
Theorem \ref{tha2})
 and  that in general a result of such kind is not true
(see Theorem \ref{tha1}).

More precisely, let $E\subset \R$ be a regular compact set. For
$\de>0$ and $z\in \C$ set $$ E_\de:=\{ z\in\Om\, :g_\Om(z)=\de\},
$$ $$ \rho_\de(z):=\mb{dist}(z,E_\de). $$ It turns out that even
for $f\in C^\al(E)$, polynomials satisfying an analogue of
(\ref{a1.2}) cannot be constructed for any $E$ under
consideration.
\begin{thm}
\label{tha1} (\cite{andntd}) There exist a regular compact set
$E_0\subset\R$ and for any $0<\al\le 1$ a function $f_\al\in
C^\al(E_0)$ such that the following assertion is false: for any
$n\in \N$ there is a polynomial $p_n\in\bP_n$ with the property:
\beq \label{a2.1} |f_\al(x)-p_n(x)|\le c\, \rho^\al_{1/n}(x),\quad
x\in E_0, \eeq where the constant $c>0$ is independent of $n$ and
$x$.
\end{thm}
 The construction of $E_0$ in Theorem \ref{tha1} uses ideas from
  Section 2. That is, let
  $$
U_0:=\{ w=\xi+i\eta:\, -\frac{\pi}{2}<\xi< \frac{\pi}{2},\,
\eta>0\}\setminus \bigcup_{k=-\infty\atop k\not= 0}^\infty J_k',
$$ where $$ J_k':=\left[
\frac{1}{k|k|},\frac{1}{k|k|}+\frac{6i}{|k|}\right]. $$ Consider
the conformal mapping $\psi_0$ of $U_0$ onto $\He$, normalized by
the boundary conditions $$ \psi_0(\infty)=\infty,\quad
\psi_0\left(\pm\frac{\pi}{2}\right)=\pm 1. $$ We extend the
inverse mapping $\phi_0:=\psi_0^{-1}$ continuously to $\ov{\He}$
(because of the symmetry of $U_0$ we have $\phi_0(0)=0$) and set
$$ J_k:=\{ x\in\R:\, \phi_0(x)\in J_k'\}, $$ $$ I_k=[x_k',x_k'']:=
\left\{\begin{array}{ll}
\psi_0(\left[-\frac{\pi}{2},-1\right]\cup\left[
1,\frac{\pi}{2}\right]), &\, k=0,\\[2ex]
\psi_0([\frac{1}{(k+1)^2},\frac{1}{k^2}]),&\, k\in\N ,\\[2ex]
\psi_0([-\frac{1}{k^2},-\frac{1}{(k-1)^2}]),&\, -k\in\N .
\end{array}\right.
$$ Then $$ E_0:=\left( \bigcup_{k=-\infty}^\infty I_k \right)
\cup\{0\} $$ satisfies the conditions of Theorem \ref{tha1}.

The analysis of the construction  above shows that $E_0$ is ``too
scarce" in a neighborhood of $0\in E_0$. Hence, to admit estimates
like (\ref{a1.2}) or (\ref{a2.1}), $E$ has to be ``thick enough"
in a neighborhood of each of its points. In order to formulate the
appropriate restrictions we need some notations.

The  set  $\R\setminus E$ consists of a finite or infinite number
of components, i.e., disjoint open intervals.

We say that $E\in\cL (\al,c),\, \al>0,c>0$, if for any bounded
component $J$ of $\R\setminus E$ the inequality
 \beq\label{4.22.2d}
\mb{dist}(J,(\R\setminus E)\setminus J)\ge c\, |J|^{1/(1+\al)}
 \eeq
 holds.

By definition, we relate a single closed interval to $\cL
(\al,c)$.

We can now state the analogue of the Nikol'skii-Timan-Dzjadyk
theorem for functions continuous on a compact subset of the real
line.
\begin{thm}
\label{tha2} (\cite{andntd})  Let the regular set $E\subset\R$
consist of a finite number of disjoint compact sets, each of which
belongs to the class  $\cL (\al,c)$ with some $\al,c>0$. Suppose
that $f\in C(E)$ and that the function $\om$ of the modulus of
continuity type satisfies (\ref{a1.1}).

Then the following conditions are equivalent:

(i) $f\in C_\om(E)$;

(ii) for any $n\in \N$ there exists a polynomial $p_n\in \bP _n$
such that $$ |f(x)-p_n(x)|\le c_8\, \om(\rho_{1/n}(x)),\quad x\in
E, $$ where the constant $c_8>0$ does not depend on $x$ and $n$.
\end{thm}
 The simplest example of $E$ satisfying the assumptions of Theorem
\ref{tha2} is the union of a finite number of disjoint closed
intervals. The compact set $$ E_\al:=\{0\}\cup
\bigcup_{n=n_\al}^\infty\left[
\frac{1}{n+1},\frac{1}{n}-\frac{1}{n^{2+\al}} \right],\quad
\al>0,\, n_\al>2^{1/\al}, $$ which obviously satisfies the
conditions of Theorem \ref{tha2}, illustrates a nontrivial
extension of (\ref{a1.2}) to compact subsets of the real line.

The proof of Theorem \ref{tha2} uses the results and ideas
concerning
 approximation of functions by complex polynomials on continua of the
special class $H^*$ introduced and discussed in the next
subsection. We outline the main steps of this proof. As usual, we
use $c_1,c_2,\cdots$ to denote positive constants that depend on
parameters inessential to the argument.

(ii)$\Rightarrow$(i). Since by our assumption  (\ref{4.22.2d}),
for any $x\in E$ and $0<\de<1$ the length of the set
$E_{x,\de}:=E\cap\{\z:\, |\z-x|\le \de\}$ satisfies
$|E_{x,\de}|\ge c_1\,\de$, the compact set $E$ is uniformly
perfect. Hence,
 (ii)$\Rightarrow$(i)
 follows from Tamrazov's inverse theorem (see \cite[p. 138]{tam}).

(i)$\Rightarrow$(ii). Let $f\in C_\om (E)$. Applying the
 procedure described, for example,
in \cite[Chapter 1]{andbeldzj} we extend $f$ continuously to $\R$
such that $$ |f(x_2)-f(x_1)|\le c_2\om(x_2-x_1),\quad
-\infty<x_1<x_2<\infty, $$ $$ f(x)=0,\quad \mb{dist}(x,E)> 3. $$
 Further we consider the Poisson integral
$$ f(z):=\frac{1}{\pi}\int\limits_{-\infty}^{\infty} \frac{y\,
f(s)\, ds}{(x-s)^2+y^2}, \quad z=x+iy\in \He , $$ which extends
$f$ harmonically to the upper half-plane $\He$.

It can be shown that for any $z_1,z_2\in \He$ such that
$|z_1|<c,\, |z_2|<c,\, |z_1-z_2|\le \de<c$, where $c>0$ is an
arbitrary (fixed) constant, we have the inequality
 \beq
\label{4.23.3} |f(z_1)-f(z_2)|\le c_3\om(\de).
 \eeq
First we consider the case $E\in\cL (\al,c)$. Without loss of
generality we assume that $E$ consists of a finite number of
components, that is, $E=\cup_{k=1}^mI_k,\, I_k=[x'_k,x_k'']$, $$
x_1'<x_1''<\cdots<x_m'<x_m''. $$ Important is that the  estimates
below do not depend on $m$. Using a linear transformation, if
necessary, we can also assume that $x_1'=-1$ and $x_m''=1$.

For $k=1,\cdots,m-1$, set $$ J_k:=[x_k'',x_{k+1}'], $$
\begin{eqnarray*}
S_k&:=&[x_k''+2it_k,x_{k+1}'+2it_k]\\ &\cup&\{ z=x+iy:\,
x_k''-t_k\le x\le x_k'',\, |z-(x_k''+it_k)|= t_k\}\\ &\cup&\{
z=x+iy:\, x_{k+1}'\le x\le x_{k+1}'+t_k,\, |z-(x_{k+1}'+it_k)|=
t_k\},
\end{eqnarray*}
where $$ t_k:=\frac{1}{3}\min(|J_k|,|I_k|,|I_{k+1}|). $$ According
to our assumption (\ref{4.22.2d}), for  sufficiently large
positive constants $c$ and $\al$ the curve $$ l:=\{ x+iy:\, -1\le
x\le1,\, y=c\, (1-|x|)^{1+\al}\} $$ satisfies the condition
 $$ \mb{dist}(S_k,l)\ge 2\,\mb{ diam}(S_k).
 $$
  We denote by $G$ the Jordan domain bounded by
$$
\partial G=L:=E\cup(\bigcup_{k=1}^{m-1}S_k)\cup l,
$$
 and let $\Om^*:=\OC\setminus\ov{G}$.

It can be proved that $ \ov{G}\in H^* $
 (for the definition of $H^*$, see
Subsection 4.3). Therefore by \cite{and88} and (\ref{4.23.3}), for
any $n\in\N$ there exists a harmonic polynomial $$
t_n(z)=\mb{Re}\sum_{j=0}^na_j\, z^j,\quad a_j\in\C, $$ (of degree
at most $n$) such that
 \beq \label{4.23.7}
|f(z)-t_n(z)|\le c_4\om(\rho_{1/n}^*(z)),\quad z\in L,
 \eeq
  where
for $z\in\C$ and $\de>0$,
 $$
\rho_\de^*(z):=\mb{dist}(z,L_\de), $$ $$ L_\de:=\{ \z:\,
|\Phi(\z)|=1+\de\}, $$ and $\Phi$ is the conformal mapping of
$\Om^*$ onto $\De:=\{ w:|w|>1\}$ normalized by the conditions $$
\Phi(\infty)=\infty,\quad \Phi'(\infty)>0. $$
  The  calculation shows that
  \beq\label{4.3p1}
\rho_\de^*(x)\le c_5\rho_\de(x),\quad x\in E,\, 0<\de<1.
 \eeq
Thus, (\ref{4.23.7}) and (\ref{4.3p1}) imply (ii).

Let now $E=\cup_{j=1}^sE_j$, where $E_j\in\cL (\al,c)$ and
$E_j\cap E_k=\emptyset$ for $j\not= k$. For each $E_j$ (consisting
of a finite number of components) we construct an auxiliary domain
$G_j$ as above and join all $G_j$ in $\He$ by smooth arcs so that
a new set $K\supset\cup_{j=1}^s\ov{G_j}$ belongs to $H^*$.

The distances from $x\in E_j$ to the $\de$-level sets of the Green
function for $E$ and for $E_j$, denoted by $\rho_{\de,E}(x)$ and
$\rho_{\de,E_j}(x)$, are equivalent, i.e., $$
\rho_{\de,E}(x)\asymp \rho_{\de,E_j}(x),\quad x\in E_j,\, 0<\de<1.
$$ The same property holds for the distances $\rho^*_{\de,E}(x)$
and $\rho^*_{\de,E_j}(x)$ from $x\in E_j$ to the $(1+\de)$-level
line for the Riemann mapping function $\Phi$ constructed for
$\OC\setminus K$ and $\OC\setminus \ov{G_j}$ respectively. That
is, $$ \rho^*_{\de,K}(x)\asymp \rho^*_{\de,\ov{G_j}}(x),\quad x\in
E_j,\, 0<\de<1. $$ Since by (\ref{4.3p1}), $$
\rho^*_{\de,\ov{G_j}}(x)\le c_6\rho_{\de,E_j}(x),\quad x\in E_j,\,
0<\de<1, $$ we have $$ \rho^*_{\de,K}(x)\le
c_7\rho_{\de,E}(x),\quad x\in E_j,\, 0<\de<1. $$ Hence, applying
(\ref{4.23.3}) and \cite{and88} we obtain (ii).

\subsection{Simultaneous approximation and interpolation of
functions on continua in the complex plane}
 Let $E\subset\C$ be a compact set with connected complement
$\Om:=\OC\setminus E$. Denote by $A(E)$ the class of all functions
continuous on $E$ and analytic in $E^0$, the interior of $E$ (the
case $E^0=\emptyset$ is not excluded).  For $f\in A(E)$ and
$n\in\N_0:=\{0,1,2,\cdots\}$ define
 $$ E_n(f,E):=\inf\limits_{p_n\in\Pn}||f-p_n||_E.$$
  By
the Mergelyan theorem (see \cite{dzj}), $$
\lim\limits_{n\to\infty}E_n(f,E)=0,\quad f\in A(E).
 $$
  The
following assertion on ``simultaneous approximation and
interpolation" quantifies a result of Walsh \cite[p. 310]{wal}:
Let $z_1,\cdots,z_N\in E$ be distinct points, $f\in A(E)$. Then
for any $n\in\N,\, n\ge N$, there exists a polynomial $p_n\in\Pn$
such that

\beq \label{4.51.1} ||f-p_n||_E\le c\, E_n(f,E),
 \eeq
  $$
p_n(z_j)=f(z_j),\quad j=1,\cdots,N,
 $$
  where $c>0$ is independent
of $n$ and $f$.

A suitable polynomial has the form $$
p_n(z)=p_n^*(z)+\sum_{j=1}^N\frac{q(z)}{q'(z_j)(z-z_j)}
(f(z_j)-p_n^*(z_j)), $$ where $$ q(z):=\prod_{j=1}^N(z-z_j), $$
and $p_n^*\in\Pn$ satisfies $$ ||f-p_n^*||_E=E_n(f,E). $$ It is
natural to ask whether it is possible to interpolate the function
$f$ as before at arbitrary prescribed points and to simultaneously
approximate it in a ``more subtle way" than in (\ref{4.51.1}). The
theorem of Gopengauz \cite{gop} about simultaneous polynomial
approximation of real functions continuous on the interval
$[-1,1]$ and their interpolation at $\pm 1$ is a useful example.
For a recent account of improvements and generalizations of this
remarkable statement (for real functions) we refer the reader to
\cite{pri}, \cite{szaver}, \cite{kilpre}.

We are going to make use of
 the D-approximation (named after Dzjadyk, who found in the
late 50s - early 60s a constructive description of H\"older
classes requiring a nonuniform scale of approximation) as a
substitute for (\ref{4.51.1}).  In \cite{and1} it is shown that
for the D-approximation to hold for a continuum $E$, it is
sufficient, and under some mild restrictions also necessary that
$E$ belongs to the class $H^*$, which can be defined as follows
(cf. \cite{and2}, \cite{and3}).

From now on we assume that $E$ is a continuum, i.e.,
$\Om:=\CC\setminus E$ is simply connected. Let diam$(E)>0$ and
$L:=\partial E$ be the boundary of $E$.

We say that $E\in H$ if any points $z,\z\in E$ can be joined by an
arc $\ga (z,\z)\subset E$ whose length $|\ga(z,\z)|$ satisfies the
condition \beq \label{4.51.3} |\ga(z,\z)| \le c\, |z-\z|,\quad
c=c(E)\ge 1. \eeq Let us compactify the domain $\Om$ by prime ends
in the sense of Carath\'eodory (see \cite{pom1}). Let
$\tilde{\Om}$ be this compactification, and let $\tilde
L:=\tilde{\Om}\setminus\Om$. Suppose that $E\in H$, then all the
prime ends $Z\in\tilde L$ are of the first kind, i.e. have
singleton impressions $|Z|=z\in L$. The circle $\{\xi:\,
|\xi-z|=r\},\, 0<r< \frac{1}{2}$diam$(E)$, contains one arc, or
finitely many arcs, dividing $\Om$ into two subdomains: an
unbounded subdomain and a bounded subdomain such that $Z$ can be
defined by a chain of cross-cuts of the bounded subdomain. Let
$\ga_Z(r)$ denote that of these arcs for which the unbounded
subdomain is as large as possible (for given $Z$ and $r$). Thus,
the arc $\ga_Z(r)$ separates the prime end $Z$ from $\infty$ (cf.
\cite{bel}, \cite{andbeldzj}).

If $0<r<R<\frac{1}{2}$diam$(E)$, then $\ga_Z(r)$ and $\ga_Z(R)$
are the sides of some quadrilateral $Q_Z(r,R)\subset \Om$ whose
other two sides are parts of the boundary $L$. Let $m_Z(r,R)$ be
the module of this quadrilateral, i.e., the module of the family
of arcs that separate the sides $\ga_Z(r)$ and $\ga_Z(R)$ in
$Q_Z(r,R)$ (see \cite{ah}, \cite{lv}).

We say that $E\in H^*$ if $E\in H$ and there exist
$c=c(E)<\frac{1}{2} $diam$(E)$ and $c_1=c_1(E)$ such that \beq
\label{4.51.4} |m_Z(|z-\z|,c)-m_{\cal Z}(|z-\z|,c)|\le c_1 \eeq
for any prime ends $Z,{\cal Z}\in\tilde L$ with their impressions
$z=|Z|,\, \z=|\cal Z|$ satisfying $|z-\z|<c$.

In particular, $H^*$ includes domains with quasiconformal boundary
(see \cite{ah}, \cite{lv}) and the classes $B_k^*$ of domains
introduced by Dzjadyk \cite{dzj}. For a more detailed
investigation of the geometric meaning of conditions
(\ref{4.51.3}) and (\ref{4.51.4}), see \cite{and3}.

We study functions defined by their $k$-th modulus of continuity
$(k\in\N)$. There are a number of different definitions of these
moduli in the complex plane (see \cite{vorpol}, \cite{tam},
\cite{dyn}, \cite{she1}). The definition by Dyn'kin \cite{dyn} is
the most convenient for our purpose.

From now on suppose $E\in H^*$.  The quantity
 $$
\om_{f,k,z,E}(\de):=E_{k-1}(f,E\cap \ov{D(z,\de)}),
 $$
  where $f\in
A(E),\, k\in\N,\, z\in E,\, \de>0$ and
 $ D(z,\de):=\{\z:\,
|\z-z|\le\de\}$
 is called the $k$-th {\it local
modulus of continuity}, and $$ \om_{f,k,E}(\de):=\sup_{z\in
E}\om_{f,k,z,E}(\de) $$ the $k$-th (global) {\it modulus of
continuity} of $f$  on $E$.

By definition, the function $w=\Phi(z)$ maps $\Om$ conformally and
univalently onto $\Delta:=\{ w:\, |w|>1\}$ and is normalized by
the conditions
 $$ \Phi(\infty)=\infty,\,\,
\Phi'(\infty)>0.
 $$
 Let
 $$ L_\de:=\{\z:\, |\Phi(\z)|=1+\de\},\quad
\de>0, $$ $$ \rho_\de(z):=\mb{dist}(z,L_\de),\quad z\in\C,\,
\de>0. $$
\begin{thm} (\cite{andprivar})
\label{4.5th1} Let $E\in H^*,\, f\in A(E),\,  k\in\N$, and let
$z_1,\cdots,z_N\in E$ be distinct points. Then for any $n\in\N,\,
n\ge N+k$, there exists a polynomial $p_n\in\Pn$ such that \beq
\label{4.51.6} |f(z)-p_n(z)|\le c_1\,
\om_{f,k,E}(\rho_{1/n}(z)),\quad z\in L,
 \eeq
  \beq \label{4.51.7}
p_n(z_j)=f(z_j),\quad j=1,\cdots,N
 \eeq with $c_1>0$ independent of
$n$.

Moreover, if $E^0\not=\emptyset$ and for $0<\de<1$,
 \beq\label{4.51.77}
 \int\limits_0^\de  \om_{f,k,E}(t)\,\frac{dt}{t}\le
c_2\, \om_{f,k,E}(\de),\quad c_2=\mb{ constant }>0,
 \eeq
  then in addition to (\ref{4.51.6})
and (\ref{4.51.7}), \beq \label{4.51.8} ||f-p_n||_K\le c_3\,
\exp(-c_4n^\alpha) \eeq for every compact set $K\subset E^0$,
where the  constants $c_3,c_4$ and $\alpha\le 1$ are independent
of $n$.
\end{thm}
 The existence  of a polynomial $p_n$ satisfying (\ref{4.51.6}) is
called a {\it D-approximation} of the function $f$  ({\it
D-property} of $E$, {\it Dzjadyk-type  theorem}). For $k>1$,
(\ref{4.51.6}) generalizes the corresponding direct theorems of
Belyi and Tamrazov \cite{beltam} ($E$ is a quasidisk) and Shevchuk
\cite{she1} ($E$ belongs to the Dzjadyk class $B_k^*$). More
detailed history can be found in these papers.

It was first noticed by Shirokov \cite{shi3} that the rate of
D-approximation may  admit significant improvement strictly inside
$E$. Saff and Totik \cite{saftot2} proved that if $L$ is an
analytic curve, then an exponential rate is achievable strictly
inside $E$, while on the boundary the approximation is
``near-best". However, even for domains with piecewise smooth
boundary without cusps (and therefore belonging to $H^*$) the
error of approximation strictly inside $E$ cannot be better than
$e^{-cn^\alpha}$ (cf. (\ref{4.51.8})) with $\alpha$ which may  be
arbitrarily small (see \cite{mai}, \cite{shitot}). In the results
from \cite{mai}, \cite{shitot} and \cite{shi2} containing
estimates of the form (\ref{4.51.8}) it is usually assumed that
$\Om$ satisfies the wedge condition. For a continuum $E\in H^*$
this condition can be violated.

 We
denote by $A^r(E),\, r\in\N$, the class of functions $f\in A(E)$
which are
 $r$-times continuously
differentiable on $E$ and set $A^0(E):=A(E)$. Keeping in mind the
Gopengauz result \cite{gop} we generalize Theorem \ref{4.5th1}  to
the case of  Hermite interpolation and simultaneous approximation
of a function $f\in A^r(E)$ and its derivatives. For simplicity we
formulate  this assertion only for the case of boundary
interpolation points and without the analog of (\ref{4.51.8}).
\begin{thm} (\cite{andprivar})
\label{4.5th2d} Let $E\in H^*,\, f\in A^r(E),\, r\in\N,\, k\in\N$,
and let $z_1,\cdots,z_N\in \partial E$ be distinct points. Then
for any $n\in\N,\, n\ge Nr+k$, there exists a polynomial
$p_n\in\Pn$ such that for $l=0,\cdots,r$,
 $$
|f^{(l)}(z)-p^{(l)}_n(z)|\le c\, \rho_{1/n}^{r-l}(z)\,
\om_{f^{(r)},k,E}(\rho_{1/n}(z)),\quad z\in L, $$
 $$
 p^{(l)}_n(z_j)=f^{(l)}(z_j),\quad j=1,\cdots,N $$
with $c$ independent of $n$.
\end{thm}
 Our next goal is to allow the number of interpolation nodes
 $N$ in Theorem \ref{4.5th1} to grow infinitely with the degree
 of the approximating polynomial $n$. To this end, we specify the
choice of points $z_1,\cdots,z_N$. In order to do it optimally
from the point of view of interpolation theory we have to require
that the discrete measure $$
\nu_N=\frac{1}{N}\sum_{j=1}^N\de_{z_j}, $$ where $\de_z$ denotes
the unit mass placed at $z$, is close to the equilibrium measure
for $E$ (for details, see \cite{saftot}). Fekete points (see
\cite{pom}, \cite{saftot}) are  natural candidates for our
purpose. Even in this case the number $N-1$ cannot be equal to the
degree $n$ of the approximating polynomial (cf. Faber's theorem
\cite{fab} claiming that for $E=[-1,1]$ there is no universal set
of nodes at which to every continuous function the Lagrange
interpolating polynomials converge in the uniform norm). However,
it was first observed by Bernstein \cite{ber} that for any
continuous function on $[-1,1]$ and any small $\ve>0$, there
exists a sequence of polynomials interpolating in the Chebyshev
nodes and uniformly convergent on $[-1,1]$, such that
$n\le(1+\ve)N$. This result was developed in several directions.
In particular, Erd\"os (see \cite{erd1} and \cite{erd2}) found a
necessary and sufficient condition on the system of nodes, for
this type of simultaneous approximation and interpolation to be
valid. We generalize the results of Bernstein and Erd\"os in the
following theorem.
\begin{thm} (\cite{andprivar})
\label{4.5th2} Let $E$ be a closed Jordan domain bounded by a
quasiconformal curve $L$. Let $f,\, r,\, k$ be as in Theorem
\ref{4.5th1} and let $z_1,\cdots,z_N\in E$ be the points of an
$N$-th Fekete point set of $E$. Then for any $\ve>0$ there exists
a polynomial $p_n\in\Pn,\, n\le (1+\ve)N$, satisfying conditions
(\ref{4.51.6}) and (\ref{4.51.7}). Moreover, if (\ref{4.51.77})
holds then in addition to (\ref{4.51.6}) and (\ref{4.51.7}) we
have (\ref{4.51.8}), and the constants $c_1,c_3,c_4$ and $\alpha$
are independent of $N$.
\end{thm}

 \subsection{Open problems} We begin with
 a question that would be a complete resolution of
the Meinardus-Varga problem on the structure of an entire function
with  geometric convergence on the positive real axis of
reciprocals of polynomials to the reciprocal of the function.

 {\bf Problem 6.} {\it Let an entire function $f$ satisfy (\ref{b1.1}). Are the
following two conditions

(i) $f$ satisfies (\ref{b1.3});

(ii) there exist $s>1$ and positive constants $c$, $\theta$, and
$r_0$ such that
 $$ ||f||_{{\cal E}_r(f,s)}\le c \,
r^\theta, \quad r\ge r_0 \:
 $$
  equivalent?}\\

Note that $ (i)\Rightarrow (ii)$ is proved in \cite{andblakov}
(cf. Theorem \ref{thb1}). The inverse implication  $
(ii)\Rightarrow (i)$ appears to be much more complicated to prove.
 One of the possibilities is to use an extension of the
classical result of Bernstein on polynomial approximation of
functions analytic in a neighborhood of a subinterval of the real
axis (see \cite{devlor}) to the case of several intervals.

 {\bf Problem 7.} {\it Let $E = \bigcup^k_{j=1} I_j$  be the union of $k$ disjoint
intervals
 $I_j = [\alpha_j,
\beta_j]\subset\R .$ Is it true that for each function $f$
satisfying (\ref{b2.5}) and (\ref{b2.6}), there exists a constant
$c > 0$ depending only on $s>1$ such that}
    \begin{equation}\label{4.6b2.7}
    E_n(f,E) \le c \, ||f||_{E^s} \, s^{-n}, \quad n \in \N
    \:?
    \end{equation}

Note that (\ref{4.6b2.7}) does not depend on the geometry of $E$.
This fact makes Problem 7 much more difficult to study compared to
the known results on polynomial approximation of functions on a
finite number of intervals (cf. \cite{fuc}).

 Our prior research (see Theorem \ref{tha2})  indicates that there exists a connection
  between the
 Nikolskii-Timan-Dzjadyk approximation theorem and the concept of
 uniformly perfect sets. We
 propose to investigate the details of this connection.

 {\bf Problem 8.}  {\it Let  $E\subset\R$ be uniformly perfect.
 Suppose that $f\in C(E)$ and that the function $\om$ of
modulus of continuity type satisfies (\ref{a1.1}). Are the
following conditions

(i) $f\in C_\om(E)$;

(ii) for any $n\in \N$ there exists a polynomial $p_n\in \bP _n$
such that $$ |f(x)-p_n(x)|\le c\, \om(\rho_{1/n}(x)),\quad x\in E,
$$ where the constant $c>0$ does not depend on $x$ and $n$

 equivalent?}

  This conjecture  is motivated by the
 connection between  uniformly perfect sets and John domains
 described in Subsection 2.2. Since the behavior of a conformal mapping of
 a John domain onto the disk is well-studied (see, for example,
 \cite{pom1}),  this can be used for constructing
 polynomials with the desired properties.

Moreover, we conjecture that uniformly perfect sets present
exactly the class of sets to which Theorem \ref{thntd} can be
generalized as in Problem 8.

 {\bf Problem 9.}  {\it Suppose that a compact set $E\subset\R$ is not uniformly perfect.
 Does it follow that for any $0<\al< 1$ there exists   a function $f_\al\in C^\al(E)$ such that the
following assertion is false: for any $n\in \N$ there is a
polynomial $p_n\in\bP_n$ with the property:
 $$
|f_\al(x)-p_n(x)|\le c\, \rho^\al_{1/n}(x),\quad x\in E,
 $$
where the constant $c>0$ is independent of $n$ and $x$?}

Next we discuss the description of classes of functions with a
given
         rate of decrease of their uniform polynomial
         approximations.
  Let $E\subset\R$ be a regular compact set and let $f\in C(E)$.
 The
  following
 fundamental problem of approximation
 theory is
  another example of the interplay between smoothness properties of a
  function $f\in C(E)$, the rate of decrease of its best polynomial
  approximations $E_n(f,E)$, and the geometry of the set $E$:
 for fixed
$\alpha>0$ describe all functions $f\in C(E)$ such that
 \begin{equation}\label{4.633.1} E_n(f,E)=O(n^{-\alpha}),\quad n\to\infty.
 \end{equation}

 For $x\in E$ and $t>0$
 let the function $\de(x,t)$ be defined by the equality
 $$\rho_{\de(x,t)}(x)=t.$$

 {\bf Problem 10.} {\it Let $E\subset\R$ be uniformly perfect,
   $f\in C(E)$, and let $0<\al<1$. Are
 the following two conditions

 (i)
 the inequality
 (\ref{4.633.1}) holds;

 (ii)  for all $x_1,x_2\in
 E$,
 $$
 |f(x_2)-f(x_1)|\le c\, \de(x_1,|x_2-x_1|)^\al,
 $$
 where $c>0$ is independent of $x_1$ and $x_2$

 equivalent?}

 So far, this conjecture is confirmed for
 $E=[a,b]$ being a closed interval  \cite[p. 265]{devlor}. In general,
 it is open.

 The proof of $(i)\Rightarrow (ii)$ is given in \cite{andup}. The inverse
 implication $(ii)\Rightarrow (i)$ needs new ideas.

  {\bf Acknowledgements.}
 This research was supported by the Alexander von Humboldt
foundation  while the author was visiting the W\"urzburg
University. The author wishes to thank the members of this
university for the pleasant mathematical atmosphere they offered
him. We also thank M. Nesterenko for helping to improve this
manuscript in many ways.

V. V. Andrievskii

 Department of Mathematical Sciences

 Kent State University

 Kent, OH 44242

e-mail: andriyev@math.kent.edu

\endddoc